\documentclass[11pt]{article}
\usepackage{latexsym}
\usepackage{fullpage,epsfig}
\usepackage{fullpage}
\usepackage{amssymb}

\usepackage{color}

\newcommand{\epsf}[1]{\epsfbox{#1}}

\newfont{\bbold}{msbm10 scaled \magstep1}
\newfont{\bbolds}{msbm7 scaled \magstep1}

\newcommand{\zs}{\mbox{\bbold Z}}
\newcommand{\zss}{\mbox{\bbolds Z}}
\newcommand{\qs}{\mbox{\bbold Q}}

\newcommand{\cs}{\mbox{\bbold C}}

\newcommand{\bx}{\bar x}
\newcommand{\by}{\bar y}

\newcommand{\F}{G}
\newcommand{\G}{F}
\newcommand{\GL}{\mathbb{L}}

\newcommand{\Ref}[1]{(\ref{#1})}
\newcommand{\beq}{\begin{equation}}
\newcommand{\eeq}{\end{equation}}
\newcommand{\gf}{generating function}
\newcommand{\gfs}{generating functions}

\def\cqfd{\par\nopagebreak\rightline{\vrule height 3pt width 5pt depth 2pt}
\medbreak}

 \newtheorem{Theorem}{Theorem}
 \newtheorem{Proposition}[Theorem]{Proposition}
\newtheorem{Corollary}[Theorem]{Corollary}
\newtheorem{Lemma}[Theorem]{Lemma}

\title{\bf{Walks in the quarter plane: Kreweras' algebraic model}}
\author{
\parbox{58mm}{
 {\sc Mireille Bousquet-M{\'e}lou\thanks{Partially
 supported by the European Community IHRP
Program, within the Research Training Network "Algebraic Combinatorics
in Europe", grant HPRN-CT-2001-00272.}}\\
{\small CNRS, LaBRI, Universit\'e Bordeaux 1 \\
351 cours de la Lib\'eration \\
33405 Talence Cedex,  France\\
{\tt mireille.bousquet@labri.fr }
}}
}

\date{}

\begin{document}
\maketitle

\begin{abstract}
We consider planar lattice walks that start from $(0,0)$, remain in
the first quadrant $i, j \ge 0$, and are made of three types of steps:
North-East, West and South. These walks are known to have remarkable
enumerative and probabilistic properties:

-- they are counted by nice numbers (Kreweras 1965),

-- the \gf \ of these numbers is algebraic (Gessel 1986),

-- the stationary distribution of the corresponding Markov chain in the
quadrant has an algebraic probability \gf \ (Flatto and Hahn 1984).

\noindent
These results are not well understood, and have been established via
complicated proofs. Here we give a uniform derivation of all of them, which
is more elementary that those previously published.
We then go further by computing the full law of the Markov chain. This
helps to delimit the border of algebraicity: the associated
probability \gf \ is no longer algebraic, unless a diagonal symmetry
holds.

Our proofs are based on the solution of certain functional equations,
which are very simple to establish. Finding purely combinatorial
proofs remains an open problem.
\end{abstract}

\section{Introduction}
\begin{figure}[bth]
\begin{center}
\epsf{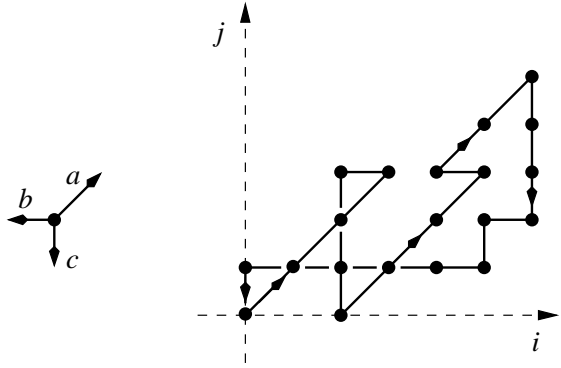}
\end{center}
\caption{Kreweras' walks in a quadrant.}
\label{figure-kreweras}
\end{figure}

Let us begin with a very simple combinatorial statement: the number of
planar lattice walks that start and end at $(0,0)$, consist of $3n$
steps that can be
North-East, South, or West, and always remain in the
nonnegative quadrant $i,j \ge 0$ is
$$
a(3n)=\frac{4^n}{(n+1)(2n+1)}{{3n}\choose n}.
$$
An example of such a walk is given in Figure~\ref{figure-kreweras}.
This result, first proved by Kreweras in 1965~\cite{kreweras}, is
rather intriguing, for at least two reasons.

 First, this simple looking formula has no simple
proof. If we consider instead the  more traditional
square lattice 
walks (consisting of North, South, East and West steps), then there
exists a nice formula too: the number of 
$2n$-step walks, starting and ending at the origin and confined in the first
quadrant, is
$$
b(2n)= \frac{1}{(2n+1)(2n+4)} {{2n+2} \choose {n+1}}^2 .$$
But the latter formula can be proved in a few lines (count first
the number of such walks having $2m$ horizontal steps, and then sum
over all values of $m$) and admits even a direct combinatorial
explanation~\cite{cori-dulucq-viennot,guy-bijections}. No similar derivation exists for the numbers $a(3n)$.

The second fact that makes the numbers $a(3n)$ intriguing is that
their  \gf , that is, the power series $A(t)=\sum_n a(3n) t^n$, is 
{\em algebraic\/}. This means that it satisfies a polynomial equation  $P(t,A(t))=0$
where $P$ is a non-trivial bivariate polynomial with rational
coefficients. For combinatorialists, objects that have an algebraic
\gf \ are 
really special: this property suggests that one {\em should\/} be able
to  factor   them into  smaller objects of
the same type, and then translate this 
factorization  into a polynomial equation (or a system of
polynomial equations) defining the  \gf. Let us take an example: it is known
that for any (finite) 
set of steps, the walks confined in the upper 
half-plane have an algebraic \gf . There is a clear combinatorial
understanding of this property: the key idea is to factor the walk at the
first time it returns to the $x$-axis. It is still an open problem to find
an explanation of this type for the algebraicity of
the series $A(t)$. Let us underline that not all
walks in the quadrant have an algebraic  \gf :  the \gf \ for 
the numbers $b(2n)$ is transcendental
(see~\cite{bousquet-petkovsek-knight} for a stronger result).

 A natural question --- at least for a computer
 scientist --- is whether the set of words on the alphabet $\{a,b,c\}$
 that naturally encode Kreweras' walks forms an {\em algebraic\/} (or:
 {\em context-free\/})  language~\cite{hopcroft}. These words contain
 as many $a$'s as $b$'s, as 
 many $a$'s as $c$'s, and each of their prefixes contains no more
 $b$'s than $a$'s, and no more $c$'s than $a$'s. Using the
{\em pumping lemma}~~\cite[Theorem 4.7]{hopcroft}, one
 can prove that this language is not 
 algebraic. Moreover, the words satisfying only
 the second 
 condition above, which encode walks ending anywhere in the quadrant, 
 do not form an algebraic language either~\cite{berstel}. However, we
 shall see that their  \gf \ is algebraic.

Then how does one prove Kreweras' formula? In his original paper,
Kreweras considered $n$-step walks in the quadrant going from $(0,0)$
to $(i,j)$. A step by step construction of these walks gives an
obvious recurrence relation for their number, denoted below
$a_{i,j}(n)$. Kreweras solved this recursion. His proof involves {\em
  guessing\/} a substantial part of the solution, and then proving
several hypergeometric identities. The latter part was then simplified
by Niederhausen~\cite{niederhausen-proceedings,niederhausen-ballot}. A
different proof, due to Gessel, also requires {\em guessing\/} the
bivariate \gf \ of walks ending on the $x$-axis, and then verifying
that it satisfies a certain functional equation~\cite{gessel-proba}. 

\medskip
On the probabilistic side, in the early 70's  Malyshev began to address
the very general problem of computing the stationary distribution of
discrete homogeneous Markov chains in the quadrant~\cite{malyshev}. Several instances of this
question actually correspond to finding the equilibrium behaviour of
 double-queue
processes~\cite{fayolle-conf,fayolle-processors,flatto-hahn,wright}.
This work culminated   in 1999 with a book  
that is entirely devoted to solving this problem in the case of
unit increments (Fayolle, Iasnogorodski and Malyshev~\cite{fayolle-livre}).
The techniques used in this book are far from elementary,  involving
sophisticated complex analysis, Riemann surfaces and boundary value
problems. Solutions are often expressed in terms of elliptic functions.
The book lists a number 
of cases in which the stationary distribution 
has a rational \gf , and mentions exactly {\em one\/} case (actually due to
Flatto and Hahn~\cite{flatto-hahn}) where this \gf \ is
algebraic. Not surprisingly, the set of increments of this
random walk is the same as in Kreweras' problem\footnote{A (partial)
  algebraicity criterion is actually given
  in~\cite[Thm.~4.3.1 and 4.3.6]{fayolle-livre}, but it is only
  illustrated by Kreweras' example.}. 

\medskip
Hence the following question: what is so special with this set of three
steps? Could one find a single argument that proves both the
algebraicity of the \gf  \ that counts  these walks and the
algebraicity of the \gf  \ for the stationary distribution of the
corresponding Markov chain?
 
This is the question we answer --- positively --- in this paper. For
both the combinatorial problem and the probabilistic one, it is very
easy to establish a functional equation defining the \gf . We solve
both equations using the same approach. The only difference is that we
are dealing with formal power series in the first problem, but with 
analytic functions in the second one.  Our solution is constructive
(we do not have to guess anything) and more elementary than the
previously published ones. In particular, we always remain in the
(small) world of algebraic functions, and do not need to introduce
elliptic functions. The key to our approach is the combination of the
{\em kernel method\/} (which is also central in~\cite{flatto-hahn}
or~\cite{fayolle-livre})  with a 
special property of the kernel of the equations we consider. Moreover,
after having solved the counting problem (Section~\ref{section-count}) and the
probabilistic one (Section~\ref{section-stationary}), we combine both
viewpoints and 
compute explicitly the full law of the Markov chain
(Section~\ref{section-law}). This actually marks the end of algebraicity: the
probability \gf \ is transcendental, unless the transition
probabilities satisfy a diagonal symmetry. Still, this \gf \ belongs
to the nice class of {\em D-finite (or: holonomic) series\/}, which
is defined below. 

\medskip
Obviously, since this paper aims at explaining why a specific set of steps
has such special properties, it cannot compete in generality with the
strength of the machinery developped in~\cite{fayolle-livre}. Still,
it is natural to ask how far our approach could be generalized.
It is actually while fighting with Kreweras' walks that I discovered
it. But it turns out that other applications of this approach were
published before I was able to complete the present paper. 
In particular, most of the results in
Section~\ref{section-count} are already reported in some conference
proceedings together with a general holonomy criterion for the
enumeration of walks in the quadrant~\cite{bousquet-versailles}.  More
recently, the same ideas were applied to certain counting problems on
permutations~\cite{bousquet-motifs}. Four new equations were thus
solved in a uniform, elementary way. Their solutions are usually
transcendental, but holonomic, and can be expressed as integrals of
algebraic (quadratic) functions. 
I have not
tried to attack other stationary distribution examples. But it is
likely that the method presented here and in~\cite{bousquet-motifs}
can be applied to solve explicitly (and in an elementary way) at least
{\em certain\/} specific examples. 

\bigskip
Let us conclude this section by giving some definitions and notation
 on formal power series.
 Given a ring $\GL$ and $k$ indeterminates
$x_1, \ldots , x_k$, we denote by
 $\GL[x_1, \ldots , x_k]$  the ring
 of polynomials 
in $x_1, \ldots , x_k$ with coefficients in $\GL$.
We denote by  $\GL[[x_1, \ldots , x_k]]$ the ring of formal power
series in the $x_i$,  that is, of formal sums
 \beq 
\sum_{n_1\ge 0, \ldots , n_k \ge 0} a(n_1, \ldots ,
 n_k)x_1^{n_1}\cdots x_k^{n_k}, 
\label{formal-series}
\eeq  
 where $ a(n_1, \ldots , n_k) \in \GL$.
 A {\em Laurent polynomial\/} in the $x_i$ is a polynomial in both the
 $x_i$ and the $\bx _i=1/x_i$. 
 A {\em Laurent series\/} in the $x_i$ is a series of the
 form~\Ref{formal-series} in which the summation runs over ${n_i\ge
 m_i}$ for all $i$, with $m_i$ in $\zs$.  
For $F\in \GL[[t]]$, we denote by $[t^n]F$ the coefficient of $t^n$ in
$F(t)$. 
 If $F$ is a formal series in $t$ whose coefficients are
Laurent series in $x$, we denote by $F^+$ the {\em positive part of $F$ in
  $x$\/}, that is,
\beq \label{positive-part}
F= \sum_{n\ge 0} t^n \sum_{i\in \zss} f_i(n) x^i \ \Rightarrow \ 
F^+= \sum_{n\ge 0} t^n \sum_{i > 0} f_i(n) x^i.
\eeq
We define similarly  the negative, nonnegative and nonpositive parts
of $F$.

Assume, from now on, that $\GL$ is a field. We denote by
 $\GL(x_1, \ldots , x_k)$ the field of rational functions in $x_1,
\ldots , x_k$ with coefficients in $\GL$.
%
A series $F$ in $\GL[[x_1, \ldots , x_k]]$ is {\em algebraic\/}
if there exists a non-trivial polynomial $P$ with coefficients in
$\GL$ such that 
$P(F,x_1, \ldots , x_k)=0.$ 
The sum and product of algebraic series is algebraic.  
%
%
The series $F$   is {\em D-finite\/} 
if the partial derivatives of $F$ span a finite 
dimensional vector space over the field $\GL(x_1, \ldots , x_k)$;
see~\cite{stanleyDF} for the one-variable case,
and~\cite{lipshitz-diag,lipshitz-df} otherwise.  In other 
words, for $1\le 
i\le k$, the series $F$ satisfies a non-trivial partial differential
equation of the form
$$\sum_{\ell=0}^{d_i}P_{\ell,i}\ \frac{\partial ^\ell
F}{\partial x_i^\ell} =0,$$
where $P_{\ell,i}$ is a polynomial in the $x_j$.
Any algebraic series is D-finite.
 The sum and product of D-finite series are D-finite.
The specializations of a D-finite series (obtained by giving 
values from $\GL$  to some of the variables) are D-finite, if
well-defined. 
 Finally, if $F$ is 
D-finite, then any {\em diagonal\/} of $F$  is also
 D-finite~\cite{lipshitz-diag} (the diagonal of $F$ in $x_1$ and $x_2$ is
 obtained by keeping only those monomials for which the exponents of
 $x_1$ and $x_2$ are equal).  We shall use the following
 consequence of the proof of this result: if
 $F(t,x) \in \GL[x,\bx ][[t]]$ is 
algebraic (with $\bx =1/x$),
then the positive part of
 $F$ in $x$ is D-finite, as well as the coefficient of $x^i$ in this
 series, for all $i$.

\section{Enumeration: the number of walks}\label{section-count}
Consider walks that start from $(0,0)$, consist of South, West and
North-East steps, and always stay in the first quadrant
(Figure~\ref{figure-kreweras}). Let 
$a_{i,j}(n)$ be the number of $n$-step walks of this type ending at $(i,j)$. We
denote by $Q(x,y;t)$ the {\em complete \gf \/} of these walks:
$$ 
Q(x,y;t):=\sum_{i,j,n \ge 0} a_{i,j}(n) x^i y^j t^n .
$$
We can construct these walks recursively, by starting from $(0,0)$   and
adding a step at each time. This gives the equation:
$$
Q(x,y;t)=1+t\left(\frac 1 x + \frac 1 y +xy\right) Q(x,y;t) -\frac t y
\ Q(x,0;t) - \frac t x \ Q(0,y;t).
$$
The first term in the right-hand side encodes the empty walk, reduced
to the point $(0,0)$. The next term shows the three possible ways one
can add a step at the end of  a walk. However, one should not add a
South step to 
a walk that ends on the $x$-axis: the third term subtracts the
contribution of this forbidden move, and the last term takes care of
the symmetric case. Equivalently,
\begin{equation}
\left( xy-t(x+y+x^2y^2)\right) Q(x,y;t) = xy-xtQ(x,0;t) -ytQ(0,y;t).
\label{qdp-kreweras1}
\end{equation}
We shall often denote $Q(x,y;t)$ by $Q(x,y)$ for short. Let us
also denote the series $xtQ(x,0;t)$ by $R(x;t)$ or even $R(x)$. Using the
symmetry of the problem in $x$ and $y$, the above equation
becomes:
\begin{equation}
\left( xy-t(x+y+x^2y^2)\right) Q(x,y) = xy-R(x) -R(y).
\label{qdp-kreweras2}
\end{equation}
 Equation~\Ref{qdp-kreweras1} is equivalent to a recurrence relation
 defining the 
numbers $a_{i,j}(n)$  inductively with respect to $n$. Hence, it
defines completely 
the series $Q(x,y;t)$. Still, the characterization  we
have in mind is of a different nature:
\begin{Theorem}[The number of walks]
\label{theorem-kreweras}
Let $W\equiv W(t)$ be the power series in $t$ defined by
$$W=t(2+W^3).$$
Then the \gf \ of Kreweras' walks ending on the $x$-axis is
$$Q(x,0;t)= \frac 1 {tx} \left( \frac 1 {2t} - \frac 1 x - 
\left( \frac 1 W -\frac 1 x \right) \sqrt{1-xW^2} \right).$$
Consequently, the length \gf \ of walks ending at $(i,0)$ is
$$[x^i] Q(x,0;t) = \frac{W^{2i+1}}{2.4^i\ t}\left( C_i
-\frac{C_{i+1}W^3}4\right),$$ 
where $C_i={{2i} \choose i}/(i+1)$ is the $i$-th Catalan number.
The Lagrange inversion formula gives the number of such walks of 
length $3n+2i$ as
$$a_{i,0}(3n+2i)=\frac{4 ^n (2i+1)}{(n+i+1)(2n+2i+1)} {{2i} \choose i}
{{3n+2i} \choose n} .$$
\end{Theorem}
The aim of this section is to derive Theorem~\ref{theorem-kreweras}
from the functional equation~\Ref{qdp-kreweras1}.  Note that the
complete \gf \ $Q(x,y)$ can the be recovered
using~\Ref{qdp-kreweras1}:
$$
Q(x,y;t)= \frac{(1/W-\bx) \sqrt{1-xW^2} +(1/W-\by ) \sqrt{1-yW^2}}
{xy-t(x+y+x^2y^2)} -\frac 1 {xyt},
$$
with $\bx =1/x$. For walks ending on the diagonal, we shall
also obtain a nice \gf : 
\begin{Theorem}[Walks ending on the diagonal]
\label{theorem-kreweras-diag}
Let $W\equiv W(t)$ be defined as above.
Then the \gf \ of Kreweras' walks ending on the diagonal, defined by
$$
Q_d(x;t):=\sum_{i,n \ge 0} a_{i,i}(n) x^i  t^n .
$$
satisfies
$$
t Q_d(x;t)=\frac {W-\bx}{\displaystyle \sqrt{1-xW(1+W^3/4)+x^2W^2/4}}+\bx.
$$
\end{Theorem}
\medskip
\noindent The expression of $Q_d$ becomes a bit simpler if we express
it in terms of the unique power series $Z\equiv Z(t)$ satisfying
$Z=1+4t^3Z^3$. Then $W=2tZ$ and
$$
t Q_d(x;t)=\frac {2tZ-\bx}{\displaystyle \sqrt{1-xtZ(1+Z)+x^2t^2Z^2}}+\bx.
$$
The last formula of Theorem~\ref{theorem-kreweras} is due to
Kreweras~\cite{kreweras}.  He also gave a closed form expression for the number
of walks containing exactly $p$ West steps, $q$ South steps, and $r$
North-East steps, that is, for walks of length $n=p+q+r$ ending at
$(i,j)=(r-p,r-q)$. This expression is a double summation, with
alternating signs. We have not found anything simpler.

\subsection{The obstinate kernel method}
\label{section-obstinate}
The kernel method is basically the only tool we have to attack
Eq.~\Ref{qdp-kreweras2}. This method has been around
since, at least, the 
70's, and
is currently the subject of a 
certain
revival (see~\cite[Ex.~2.2.1.4 and 2.2.1.11]{knuth} and
\cite{fayolle-conf} for early uses of the method,
and~\cite{hexacephale,banderier-flajolet,bousquet-petkovsek-recurrences}
for more recent combinatorial applications). 
It consists in coupling the
variables $x$ and $y$ so as to cancel the {\em kernel\/}  $K(x,y)=
xy-t(x+y+x^2y^2)$ (which is the coefficient of $Q(x,y)$
in~\Ref{qdp-kreweras2}). 
This should give the ``missing'' information about
the series $R(x)$.

As a polynomial in $y$, this kernel has two roots
$$
\begin{array}{lclllllll}
Y_0(x)&=&\displaystyle \frac{1-t\bx -\sqrt{(1-t\bx )^2-4t^2x}}{2tx}&=&
 & &t+\bx t^2 + O(t^3) , \\
\\
 Y_1(x)&=&\displaystyle \frac{1-t\bx +\sqrt{(1-t\bx )^2-4t^2x}}{2tx}
&=&\displaystyle\frac \bx t -\bx ^2 &-&t-\bx t^2 + O(t^3) .
\end{array}$$
The elementary symmetric functions of the $Y_i$ are
\begin{equation}
Y_0+Y_1= \frac \bx {t} - \bx ^2 \quad \hbox{and}
\quad Y_0Y_1 = \bx .
\label{symmetric-functions}
\end{equation}
The fact that they are polynomials in $\bx =1/x$ will play a very
important role below.

Only the first root can be substituted for $y$
in~\Ref{qdp-kreweras2} (the term $Q(x,Y_1;t)$ is not a well-defined
power series in $t$, because of the negative power of $t$ that occurs
in $Y_1$). We thus obtain a functional equation for
$R(x)$:
\begin{equation}
R(x)+R(Y_0)=xY_0.
\label{kernel0}
\end{equation}
It is not hard to see that this equation --- once restated in terms of
$Q(x,0)$ --- defines uniquely $Q(x,0;t)$ as a
formal power series in $t$ with polynomial coefficients in
$x$. Equation~\Ref{kernel0} is the standard result of the kernel method.

Still, we want to apply here the {\em obstinate\/} kernel method. That
is, we shall not content ourselves with Eq.~\Ref{kernel0}, but we
shall go on producing pairs $(X,Y)$ that cancel the kernel and use the
information they provide on the series $R(x)$.
This obstinacy was inspired by the book~\cite{fayolle-livre} by Fayolle,
Iasnogorodski and Malyshev, and more precisely by Section~2.4 of this
book, where one possible way to obtain such pairs is described (even
though the analytic context is different). We give here
an alternative construction.

\begin{figure}[hbt]
\begin{center}
\input{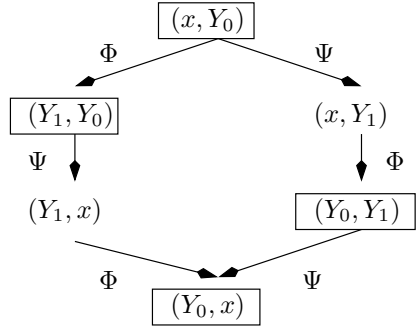}
\end{center}
\caption{The orbit of  $(x,Y_0)$ under the action of  $\Phi$ and
  $\Psi$. The framed pairs can be  substituted for $(x,y)$ 
in the functional equation.}
\label{diagram}
\hrule
\end{figure}

Let $(X,Y)\not = (0,0)$ be a pair of Laurent series in $t$ with
coefficients in some field  such that $K(X,Y)=0$. Recall that, as a
function of $y$, the polynomial $K(x,y)$ is quadratic. Thus let  $Y'$
be {\em the other solution\/}  of the equation $K(X,y)=0$. We define
the function $\Psi$ by  $\Psi(X,Y)= (X,Y')$. For instance, if $(X,Y)$
is the pair $(x,Y_0)$, then $\Psi(X,Y)=(x,Y_1)$. Similarly, we define
$\Phi(X,Y)= (X',Y)$, where $X'$ is  the other solution  of $K(x,Y)=0$.
Note that $\Phi$ and $\Psi$ are involutions and that, in view
of~\Ref{symmetric-functions}, $X'=Y'=(XY)^{-1}$. In particular,
$\Phi(x,Y_0)=(Y_1,Y_0)$.  Let us examine the iterated action of $\Phi$
and $\Psi$ on the pair $(x,Y_0)$: We obtain the diagram of
Figure~\ref{diagram}. 

\label{page-diagram}
All these pairs of power series cancel the kernel. We have already
seen that the pair $(x,Y_0)$ can be substituted for $(x,y)$ in
Eq.~\Ref{qdp-kreweras2}. This is also true, but less obvious, for the
pair $(Y_0,Y_1)$: indeed, if we write
\begin{eqnarray*}
Q(x,y;t)&=& \sum_{k \ge \max(\ell, m) } t^{k+\ell +m} x^{k-\ell} y^{k-m}
a_{k-\ell, k-m}(k+\ell+m),\\
&=&
\sum_{k \ge \max(\ell, m)} t^{k +2m} (xy)^{k-\ell} (ty)^{\ell-m}
a_{k-\ell, k-m}(k+\ell+m),
\end{eqnarray*}
and note that $Y_0Y_1=\bx$, while $tY_1=\bx + O(t)$, then we see that
$Q(Y_0,Y_1;t)$ is a well-defined power series in $t$, with
coefficients in $\qs[x,\bx]$. The same argument shows that
$R(Y_1)=tY_1Q(0,Y_1;t)$ is also well-defined.
 Thus the two pairs than can be substituted for $(x,y)$ in the
 functional equation give us
 {\em two\/} equations for the unknown series $R(x)$:
\beq
\label{kernel2.1}
\left\{
\begin{array}{lll}
R(x)+R(Y_0)&=&xY_0, \\
R(Y_0)+R(Y_1)&=&Y_0Y_1=\bar x. 
\end{array}
\right.
\eeq

\medskip
\noindent{\bf Remark.} Let $p,q,r$ be three nonnegative numbers such that
$p+q+r=1$. Take $x= (pr)^{1/3}q^{-2/3}$,  $y= (qr)^{1/3}p^{-2/3}$, and
$t=(pqr)^{1/3}$. Then $K(x,y;t)=0$, so that $R(x)+R(y)=xy$. This
equation can be given a probabilistic interpretation by considering
random walks that make a North-East step with (small) probability $r$
and a West (resp. South) step  with probability $p$ (resp.~$q$). This
probabilistic argument, and the equation it implies, is the starting
point in Gessel's solution of Kreweras 
problem~\cite[Eq.~(21)]{gessel-proba}.

\subsection{Symmetric functions of $Y_0$ and $Y_1$}
\label{section-symmetric}
After the kernel method, the next tool in our approach is the extraction
of the positive part of power series, defined by~\Ref{positive-part}.
This is where the values of the symmetric functions of $Y_0$ and $Y_1$
become crucial: the fact that they only involve negative powers of $x$
(see~\Ref{symmetric-functions}) will simplify the extraction of the positive
part of certain equations.
\begin{Lemma}
\label{lemma-symmetric}
Let $F(u,v;t) $ be a Laurent series in $t$ with coefficients in
$\cs[u,v]$, symmetric in $u$ and $v$. That is,
$F(u,v;t)=F(v,u;t)$. Then the series 
$F(Y_0,Y_1;t)$, if well-defined, is a Laurent series in $t$ with
polynomial coefficients in $\bx$. Moreover, the constant term of this
series, taken with respect to $\bx$, is $F(0,0;t)$.
\end{Lemma}
\noindent {\bf  Proof.} By linearity, it suffices to check this when
$F$ is simply a symmetric polynomial in $u$ and $v$.  But then it is a
polynomial in $u+v$ and $uv$ with complex coefficients. The result
follows, thanks to~\Ref{symmetric-functions}.
\cqfd

We now want to form a symmetric function of $Y_0$ and $Y_1$, starting
from the equations~(\ref{kernel2.1}). The first one
reads
$$R(Y_0)-xY_0=-R(x).$$
By combining both equations, we  obtain the companion expression:
$$R(Y_1)-xY_1=R(x)+ 2\bx -1/t.$$
Taking the difference\footnote{An alternative derivation of Kreweras'
result, obtained by considering the product
$(R(Y_0)-xY_0)(R(Y_1)-xY_1)$, is presented
in~\cite{bousquet-versailles}.} and dividing by $Y_0-Y_1$ gives  
\beq
\frac{R(Y_0)-R(Y_1)}{Y_0-Y_1}-x = tx \ \frac{2R(x)+2\bx
-1/t}{\sqrt{\Delta(x)}},
\label{divided-diff}
\eeq
where $\Delta(x)=(1-t\bx)^2-4t^2x$ is the discriminant that occurs in
both $Y_0$ and $Y_1$. 

As a Laurent polynomial in $x$, $\Delta(x)$ has three roots. Two of
them, say $X_0$ and $X_1$,  are formal
power series in $\sqrt t$; the other is a Laurent series in $t$ (for
generalities on the roots of a 
polynomial over $\cs(t)$, see~\cite[Chapter 6]{stanley-vol2}). The
coefficients of these series can be computed inductively:
\begin{eqnarray*}
X_0 & = &  t+2t^2\sqrt t+6t^4+21t^5\sqrt t+80t^7+\frac{1287}{4}t^8 \sqrt t+ \cdots \\
X_1 & = &   t-2t^2\sqrt t+6t^4-21t^5\sqrt t+80t^7-\frac{1287}{4}t^8
\sqrt t + \cdots \\
X_2 & =&  \frac 1 {4t^2}
-2t-12t^4-160t^7-2688t^{10}-50688t^{13} + \cdots 
\end{eqnarray*}
Hence $\Delta(x)$ factors as
$$
\Delta(x)=\Delta_0 \Delta_+(x) \Delta_-(\bx)
$$
with
$$
\Delta_0= 4t^2 X_2, \quad \Delta_+(x)= 1 -x/X_2, \quad \Delta_-(\bx)=
(1-\bx X_0)(1-\bx X_1).
$$
Note that $\Delta_0 ,\Delta_+(x)$ and $ \Delta_-(\bx)$ are  power
series in $t$ with constant term $1$. Moreover,  $\Delta_0$ has its
coefficients in $\qs$, while $\Delta_+(x)$ has its coefficients in
$\qs[x]$, and $ \Delta_-(\bx)$ has its  coefficients in
$\qs[\bx]$. This is an instance of the ``canonical factorization'' of
power series 
of $\qs[x, \bx] [[t]]$, which has already proved useful in several
path enumeration
problems~\cite{gessel-factorization,bousquet-slit,bousquet-schaeffer-slit}.
Going back to~\Ref{divided-diff}, 
and multiplying through by $\sqrt{\Delta_-(\bx)}$, one obtains
$$
\sqrt{\Delta_-(\bx)} \left( \frac{R(Y_0)-R(Y_1)}{Y_0-Y_1}-x\right) =
t \ \frac{2xR(x)+2-x/t}{\sqrt{\Delta_0 \Delta_+(x)}}.
$$
Both sides of this identity are power series in $t$ with coefficients
in $\qs[x,\bx]$. But the right-hand side only contains nonnegative
powers of $x$, while the left-hand side, except for a term $-x$, only
contains nonpositive powers of $x$ (in view of
Lemma~\ref{lemma-symmetric}). Extracting the positive part of the
above equation thus gives
$$
-x = \frac t {\sqrt{\Delta_0}} \left(
\frac{2xR(x)+2-x/t}{\sqrt{\Delta_+(x)}}-2\right).
$$
The expression of $Q(x,0)$  announced in
Theorem~\ref{theorem-kreweras} follows, given that 
 $X_2=1/W^2$ and $R(x)=xtQ(x,0)$. 
The expansion of $Q(x,0)$ in $x$ is straightforward, using
$1-\sqrt{1-4t}=2t \sum_{n\ge 0} C_n t^n$. The value of $a_{i,0}(3n+2i)$
follows using the Lagrange inversion
formula~\cite[p.~38]{stanley-vol2}.\cqfd

\subsection{The algebraic kernel method}
\label{section-algebraic-kernel}
We present in this section another proof of
Theorem~\ref{theorem-kreweras} based on a variation of the kernel
method. This variation does not require to cancel the kernel, but,
instead, builds on one of its algebraic properties. This variant has
some drawbacks --- since the kernel is not zero, we are handling
bigger equations --- but it also has some advantages. In particular, we 
obtain at some point an equation that is the counterpart
of~\Ref{divided-diff}, but in which it is obvious that the left-hand
side is nonpositive in $x$. This will be helpful in the next section,
where we handle analytic functions rather than power series. Finally,
this variant of the kernel method provides a proof of
Theorem~\ref{theorem-kreweras-diag}.

 Let us return to the original equation~\Ref{qdp-kreweras2}, or,
equivalently, to
$$
xyK_r(x,y)Q(x,y)=xy-R(x) -R(y),
$$
where $K_r(x,y)=1-t(\bx +\by +xy)$ is the {\em rational version\/} of the
kernel $K$. The fact that the 
diagram of Figure~\ref{diagram} is nice actually
stems from an invariance property of $K_r$:
$$
K_r(x,y)=K_r(\bx \by ,y)=K_r(x,\bx \by)\equiv K_r.
$$
Applying iteratively the (involutive) transformations $\Phi: (x,y) \mapsto (\bx \by
,y)$ and  $\Psi: (x,y) \mapsto (x,\bx \by
)$ gives the following set of pairs, on which   $K_r$ takes the same value:

\begin{figure}[hbt]
\begin{center}
\input{diagram-formel-kreweras.pstex_t}
\end{center}
\label{diagram-formel-kreweras}
\end{figure}

\hskip -5mm
\noindent Note that the diagram of Figure~\ref{diagram} is the
specialization of the above one to the case $y=Y_0$. 

Now, {\em all\/} pairs of the above diagram can be substituted for
$(x,y)$ in 
the functional equation: the resulting series are power series in $t$
with  coefficients in $\qs[x,\bx, y, \by]$. This gives no less than
{\em three\/} equations:
$$
\begin{array}{rcccccccccccccccll}
xy\ K_r&Q(x,y)&=&xy&-&R(x) &-&R(y), \\
\bx \ K_r&Q(\bx \by ,y)&=&\bx &-&R(\bx \by ) &-&R(y), \\
\by \ K_r&Q(x,\bx \by)&=& \by&-&R(x) &-&R(\bx \by).
\end{array}
$$
We sum the first and third equations, and subtract the second one, so
as to keep $R(x)$ as the only unknown function on the right-hand side:
$$
K_r\Big( xyQ(x,y)-\bx Q(\bx \by ,y)+\by Q(x,\bx \by)\Big)= xy -\bx +\by
-2R(x) = \frac {1-K_r} t -2\bx  -2R(x).
$$
Equivalently,
\beq
xyQ(x,y)-\bx Q(\bx \by ,y)+\by Q(x,\bx \by)+\frac 1 t
= \frac 1 {K_r} \left( \frac 1 t -2\bx  -2R(x)\right).
\label{noyau-formel1}
\eeq
 The kernel $K(x,y)$ factors as
$-tx^2(y-Y_0)(y-Y_1)$. Converting $1/K$ into partial fractions of
$y$ yields the following expression for the reciprocal of the
 (rational) kernel $K_r$:
$$
\frac 1 {K_r} = \frac 1 {\sqrt{\Delta(x)}} \left( \frac 1 {1-\by Y_0} +
\frac 1 {1-y /Y_1}-1\right)
=  \frac 1 {\sqrt{\Delta(x)}} \left( \sum_{n\ge 0} \by ^n Y_0^n +
 \sum_{n\ge 1} y ^n Y_1^{-n}\right).
$$
Note that this expansion is valid in the set of formal power series in
$t$
with coefficients in $\qs[x,\bx,y,\by]$.
Let us extract in~\Ref{noyau-formel1} the constant term in $y$: the
series $xyQ(x,y)$ and $\by Q(x,\bx \by)$ do not contribute, and we obtain
$$
-\bx Q_d(\bx)+\frac 1 t= \frac { 1/ t
-2\bx  -2R(x)} {\sqrt{\Delta(x)}}
$$
where the series $Q_d$ is the diagonal of $Q(x,y)$, and counts walks
ending on the diagonal.
The above equation should be compared to~\Ref{divided-diff}: basically,
both equations are equivalent, but their negative parts on the
left-hand side are written in
 two different ways. We now proceed as above, using the
canonical factorization of $\Delta(x)$, which gives
$$
\sqrt{\Delta_-(\bx)} \left(\frac 1 t -\bx Q_d(\bx)\right) =
 \frac { 1/ t -2\bx  -2R(x)} {\sqrt{\Delta_0 \Delta_+(x)}}.
$$
Extracting the nonnegative part gives, as before, the value of
$R(x)$, and Theorem~\ref{theorem-kreweras}. Extracting the negative part gives 
$$
\sqrt{\Delta_-(\bx)} \left(\frac 1 t -\bx Q_d(\bx)\right) -\frac 1 t=
-\frac {2\bx}{\sqrt{\Delta_0}}.
$$
Recall that $\Delta_0=4t^2X_2=4t^2/W^2$ and $\Delta_-(\bx)=(1-\bx
X_0)(1-\bx X_1)$, where $X_0$ and 
$X_1$ are the two ``small'' roots of $\Delta(x)$. We can express their
elementary symmetric functions in terms of the third root,
$X_2=1/W^2$. This gives 
$$
 \Delta_-(\bx) =1-\bx W(1+W^3/4)+\bx^2W^2/4,
$$
and this provides the expression of $Q_d(x)$
given in Theorem~\ref{theorem-kreweras-diag}.
\cqfd

\section{Probability: a Markov chain and its stationary distribution}
\label{section-stationary}

\begin{figure}[htb]
\begin{center}
\input{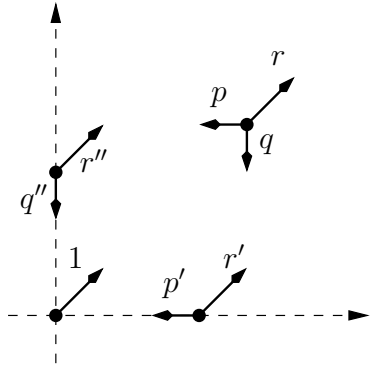}
\end{center}
\caption{The transition probabilities.}
\label{figure-transitions}
\end{figure}

We consider a Markov chain on the quadrant, whose  transition
probabilities $T(i,j;k,l)$  are schematized in
Figure~\ref{figure-transitions}. More precisely, for $i>0$ and $j>0$, the
probability of going from $(i,j)$ to $(k,\ell)$ is
$$T(i,j;k,\ell)=\left\{ 
\begin{array}{lll}
p & \hbox{if } k=i-1 &\hbox{ and } \ell =j \\
q & \hbox{if } k=i& \hbox{ and } \ell =j-1 \\
r & \hbox{if } k=i+1 &\hbox{ and } \ell =j+1 
\end{array}
\right.
$$
where $p,q,r$ are three positive real numbers summing to $1$. When
the point $(i,j)$ lies on the border of the quadrant, the transition
probabilities are modified as follows: for $i>0$, 
$$T(i,0;k,\ell)=\left\{ 
\begin{array}{lll}
p' & \hbox{if } k=i-1 &\hbox{ and } \ell =0 \\
r' & \hbox{if } k=i+1 &\hbox{ and } \ell =1 
\end{array}
\right.
$$
and for $j>0$,
$$T(0,j;k,\ell)=\left\{ 
\begin{array}{lll}
q'' & \hbox{if } k=0 &\hbox{ and } \ell =j-1 \\
r'' & \hbox{if } k=1 &\hbox{ and } \ell =j+1 
\end{array}
\right.
$$
where $p',r', q'',r''$ are positive numbers such that
$p'+r'=q''+r''=1$. Finally, we take $T(0,0;1,1)=1$. Note that this
chain is irreducible (all states communicate) and has period $3$.

A probability distribution $(p_{i,j})_{i,j \ge 0}$
is {\em stationary\/} for the above transition if  for all $k, \ell \ge 0$,
$$
p_{k,\ell} = \sum_{i,j} p_{i,j} T(i,j;k, \ell).
$$
Our objective is to find the stationary distribution of the above
transition, when it exists. 
It is customary to encode  a distribution by its probability \gf \ 
$$
\Pi(x,y)= \sum_{i,j \ge 0} p_{i,j} x^i y^j,
$$
but is is more convenient here to split $\Pi(x,y)$ into four parts:
firstly, $p_{0,0}$, and then the  three following  \gfs :
$$
P(x,y)= \sum_{i,j \ge 1} p_{i,j} x^i y^j, \quad 
P_1(x)= \sum_{i \ge 1} p_{i,0} x^i, \quad 
P_2(y)= \sum_{j \ge 1} p_{0,j}  y^j. 
$$
Then the distribution $(p_{i,j})_{i,j \ge 0}$ is stationary if and only if
\beq
(1-p\bx -q\by -rxy) P(x,y) +(1-p'\bx -r' xy) P_1(x) 
+(1-q''\by -r''xy) P_2(y)+(1-xy) p_{0,0}=0.
\label{eq-stat}
\eeq
Note that the numbers $p_{i,j}$ have to sum to $1$: hence the above
series are absolutely convergent for $|x|\le 1$ and $|y|\le 1$, and
define analytic functions for  $|x|< 1$, $|y|<1$. Moreover, 
\beq
p_{0,0}+ P_1(1)+P_2(1)+P(1,1)=1.
\label{normalization}
\eeq

\subsection{The main results}
The stationary distribution of this Markov chain was computed
in~\cite{fayolle-livre} in the case where the transition probabilities are
related by
\beq
\frac p r = \frac {p'}{r'} :=\alpha \quad \hbox{and} \quad 
\frac q r = \frac {q''}{r''} :=\beta.
\label{conditions}
\eeq
Equivalently,
$$
p'= \frac p{p+r}, \quad r'= \frac r{p+r}, \quad 
q''= \frac q{q+r}, \quad r''= \frac r{q+r}.
$$
It is known that this chain has a stationary distribution if and only
if $r <\min(p,q)$ (see~\cite{malyshev1,fayolle-livre-rouge} for
general results on Markov chains in the quadrant). It will be shown 
in Lemma~\ref{lemma-convergence} that the condition is necessary, and
it will follow from the 
results of Section~\ref{section-law}, where we compute the law of the
chain, that it is sufficient. For the moment, we rely on the general results
of~\cite{malyshev1,fayolle-livre-rouge}. 

 Under the conditions~\Ref{conditions}, which we assume to hold in
 this section,  we have
\begin{eqnarray*}
1-p'\bx -r'xy&=& \frac{r'} r \left( 1-p\bx -q\by -rxy +q(\by -1)\right),\\
1-q''\by -r''xy&=& \frac{r''} r \left( 1-p\bx -q\by -rxy +p(\bx
-1)\right),\\
1-xy&=& \frac1 r \left( 1-p\bx -q\by -rxy +p(\bx
-1) +q(\by-1) \right)
\end{eqnarray*}
so that the functional equation~\Ref{eq-stat} can be nicely rewritten as
\beq 
( 1-p\bx -q\by -rxy) Q(x,y) = q(1-\by) Q(x,0) + p(1-\bx) Q(0,y)
\label{eq-flatto}
\eeq
with
\beq
Q(x,y) = p_{0,0}+ r'P_1(x)+r''P_2(y)+rP(x,y).
\label{linkPQ}
\eeq
This equation was first met by Flatto and Hahn~\cite{flatto-hahn}
in their study of a system of two parallel queues with two
demands (with continuous time).
%
%
They solved this equation using non-trivial complex analysis,
multivalued analytic functions, and a parametrization of the kernel by
elliptic functions --- to end up with an {\em algebraic\/} solution
$Q(x,y)$. We shall rederive their result in a more elementary way,
and state it in a more symmetric fashion.

\begin{Theorem}[Solution of Flatto and Hahn's equation]
\label{theorem-stationary}
Assume $r < \min(p,q)$.  There exists, up to a multiplicative constant, a
unique solution of~{\em \Ref{eq-flatto}} that is analytic in $\{|x|,
|y|<1\}$ and whose series expansion converges for $|x|, |y|\le
1$. This solution satisfies:
$$
Q(x,0)=\frac{Q(0,0)}{(1-qx/p)(1-rx/p)}
\left( \left(1-\frac x {pw}\right) \sqrt{1-xqrw^2} -\frac{qx}{p}
 \left(1-\frac 1 {qw}\right) \sqrt{1-prw^2} \right) ,
$$
$$
Q(0,y) =\frac{Q(0,0)}{(1-py/q)(1-ry/q)}
\left( \left(1-\frac y {qw}\right) \sqrt{1-yprw^2} -\frac{py}{q}
 \left(1-\frac 1 {pw}\right) \sqrt{1-qrw^2} \right),
$$
where $w$ is the smallest positive solution of
$w=2+pqrw^3.$ 
The complete \gf \ $Q(x,y)$ can be obtained using~{\em\Ref{eq-flatto}}. 
When $p=q$, then $w=1/p$, and the  above expressions simplify to
$$
Q(x,0)=Q(0,x)= \frac{Q(0,0)}{\sqrt{1-rx/p}}.
$$
If $r\ge \min(p,q)$, no solution of~{\em\Ref{eq-flatto}} converges on $|x|,
|y|\le 1$. 
\end{Theorem}
Observe that exchanging $p$ and $q$, and $x$ and $y$, leaves the
solution unchanged, in conformity with the diagonal symmetry of the  
model. 
From the algebraic equation defining $w$, it is not difficult to see
that
$$
(pw-1)^2(1-qrw^2)= (qw-1)^2(1-prw^2)=(rw-1)^2(1-pqw^2).
$$
Moreover, an elementary study  of the function $f(z)=pqrz^3-z+2$ gives
bounds for $w$:
\beq
\frac 1 {M} \le w \le \frac 1{\sqrt {pq}} \le \frac 1
{m} < \frac 1 r, 
\label{ineqwpqr}
\eeq
where $m= \min(p,q)$ and $M=\max(p,q)$, with equalities holding if and
only if $p=q$. Hence
\beq
0\le (Mw-1)\sqrt{1-mrw^2}= -(mw-1)\sqrt{1-Mrw^2}=-(rw-1)\sqrt{1-pqw^2},
\label{double2}
\eeq
and this allows us to rewrite the second part of the expressions of
$Q(x,0)$ and $Q(0,x)$ in various ways. This will be useful in
Section~\ref{section-comments}, where we make further comments  on 
 this solution, and relate it to
Flatto and Hahn's formulation. We shall
prove Theorem~\ref{theorem-stationary} in
Section~\ref{section:proof}. For the 
moment, let us derive from it the stationary distribution of the
Markov chain of Figure~\ref{figure-transitions}. This result is
actually not given explicitly in~\cite{fayolle-livre}.
\begin{Corollary}[The stationary distribution]
\label{coro-stationary}
Assume $r< \min(p,q)$. Let $w$ be the smallest positive solution of
$w=2+pqrw^3.$  The Markov chain schematized in
Figure~{\em \ref{figure-transitions}}, with the additional condition {\em
\Ref{conditions}}, has a unique stationary distribution $(p_{i,j})$,
given by 
\begin{eqnarray*}
p_{0,0}&=&\displaystyle \frac{w(p-r)(q-r)|p-q|}{6pq(1-rw)\sqrt{1-pqw^2}},\\
\displaystyle \sum_{i >0} p_{i,0} x^i &=&\displaystyle \frac 1 {r'} \left( Q(x,0)-Q(0,0) \right),\\
\displaystyle \sum_{j >0} p_{0,j} y^j& =&\displaystyle \frac 1 {r''} \left( Q(0,y)-Q(0,0)\right), \\
\displaystyle \sum_{i >0,j>0} p_{i,j} x^i y^j&=&\displaystyle \frac 1 r \left (
Q(x,y)-Q(x,0)-Q(0,y)+Q(0,0) \right) 
\end{eqnarray*}
where $Q(x,y)$ is the function of Theorem~{\em \ref{theorem-stationary}},
taken with  $Q(0,0)=p_{0,0}$.  
When $p=q$, then  the expression of $p_{0,0}$ should be taken to be
$$
p_{0,0}= \frac 1 3 (1-r/p)^{3/2}.
$$
If $r \ge \min(p,q)$, then the Markov chain has no stationary distribution.
\end{Corollary}
{\bf Proof.} The stationary distribution is related to the series
$Q(x,y)$ satisfying~\Ref{eq-flatto}  by~\Ref{linkPQ}, so that the
expressions of the three series above are obvious. We only 
 have to determine which value of $Q(0,0)$  guarantees the
normalizing condition~\Ref{normalization}. This condition reads
$$
Q(0,0) + \frac 1 {r'} (Q(1,0)-Q(0,0)) + \frac 1 {r''} (Q(0,1)-Q(0,0))
+\frac 1 r (Q(1,1)-Q(1,0)-Q(0,1)+Q(0,0)) =1,
$$
that is,
\beq 
Q(1,1)-qQ(1,0)-pQ(0,1)=r.
\label{normalizingQ}
\eeq
When $y=1$, a factor $(1-\bx)$ comes out of~\Ref{eq-flatto}, leaving
\beq
(p-rx)Q(x,1)=pQ(0,1).
\label{y=1}
\eeq
Setting $x=1$  gives
$$
Q(1,1)= \frac p{p-r} Q(0,1),
$$
and of course, a symmetric argument yields
$$
Q(1,1)= \frac q{q-r} Q(1,0).
$$
Hence the normalizing condition~\Ref{normalizingQ} reads
\beq
Q(1,0)= \frac{q-r}{3q}.
\label{normalizingQ10}
\eeq
Now, from the expression of $Q(x,0)$ given in
Theorem~\ref{theorem-stationary}, we obtain, if $p\not = q$,
$$
Q(1,0)=\frac{pQ(0,0)}{w(p-q)(p-r)}
\left( \left( {pw}-1\right) \sqrt{1-qrw^2} -
 \left(qw-1\right) \sqrt{1-prw^2} \right) .
$$
Using~\Ref{double2}, the above expression for $Q(1,0)$ can be rewritten as
$$
Q(1,0)= \frac{2pQ(0,0) (1-rw)\sqrt{1-pqw^2}}{w|p-q|(p-r)},
$$
and the condition~\Ref{normalizingQ10} gives the value of
$Q(0,0)=p_{0,0}$.

When $p=q$, the simplified expression of $Q(x,0)$, given in
Theorem~\ref{theorem-stationary}, gives $Q(1,0)= Q(0,0)/\sqrt{1-r/p}$, 
and the result follows.
\cqfd

\subsection{Proof of Theorem~\ref{theorem-stationary}}
\label{section:proof}
Our solution of Eq.~\Ref{eq-flatto} follows the same idea as 
 Section~\ref{section-algebraic-kernel}: we shall exploit
an invariance property of the kernel. 
However, we do not have the length variable $t$ any more,
which means that we are no longer in a power series context, but rather in the
world of functions of two  complex variables $x$ and
$y$. Consequently, certain operations that were performed formally in
Section~\ref{section-algebraic-kernel} (e.g., the extraction of
coefficients) now need to be justified analytically. The analytic
lemmas we need are gathered in Section~\ref{analytic-lemmas} below. In
Section~\ref{section-simplification}, our main functional
equation~\Ref{eq-flatto} is transformed 
into {\em two\/} functional equations defining two functions $U(x,y)$
and $F(x,y)$, which are respectively symmetric and anti-symmetric in
$x$ and $y$. Section~\ref{section-algebraic-analytic} is the heart of
the proof: there we solve 
these two equations, using the algebraic kernel method. The reader may
skip directly to the latter section in order to recognize the logic of the
kernel method, transferred to an analytic context.

\subsubsection{Preliminary results}\label{analytic-lemmas}
Our first lemma tells us that the domain of convergence of  $Q(x,y)$
is actually larger than the unit polydisc $|x|\le 1$, $|y|\le 1$.
\begin{Lemma}
\label{lemma-convergence}
Let $Q(x,y)$ be a power series solution of~{\em \Ref{eq-flatto}} that has
nonnegative coefficients and converges for 
$|x|\le 1$ and $|y| \le 1$. Then this series also converges in the
following domain:
$$
\{ |x| <p/r, |y|\le 1\} \cup \{|x|\le 1, |y| <q/r\}.
$$
Moreover,
$$
Q(x,1)= \frac {Q(0,1)} {1-rx/p} \quad \hbox{and} \quad
Q(1,y)= \frac {Q(1,0)} {1-ry/q}.
$$
In particular, such a solution $Q(x,y)$ can only exist if $r<\min
(p,q)$: if this inequality does not hold, the Markov chain has no stationary
distribution. 
\end{Lemma}
{\bf Proof.} For $|x|\le 1$, the expression of
$Q(x,1)$ follows  from~\Ref{y=1}. The expression of $Q(1,y)$ is of
course symmetric. These two series must converge when $x=1$ and $y=1$:
this forces $r$ to be smaller than $p$ and $q$. Given the
nonnegativity of the coefficients, the values of $Q(x,1)$ and $Q(1,y)$
imply the convergence of $Q(x,y)$ in the desired domain.
\cqfd

\medskip
The extraction of coefficients will be based on the following
result~\cite[Chap.~10, Ex.~25]{rudin}. 
\begin{Proposition}\label{annulus}
Let $f(z)$ be an analytic function in the annulus ${\cal A}=\{r
<|z|<R\}.$ There exists a unique bi-infinite sequence $(a_n)_{n\in
  \zss}$ such that for all $z \in {\cal A}$,
$$
f(z)=\sum_{n\in \zss} a_n z^n.
$$
Moreover, the convergence is absolute. In other words, $f$ is the sum
of a function analytic in the disk $\{|z|<R\}$ and a function analytic for
$|z|>r$.
\end{Proposition}

\begin{figure}[ht]
\begin{center}
\epsfxsize=63mm \epsfysize=6cm
\epsf{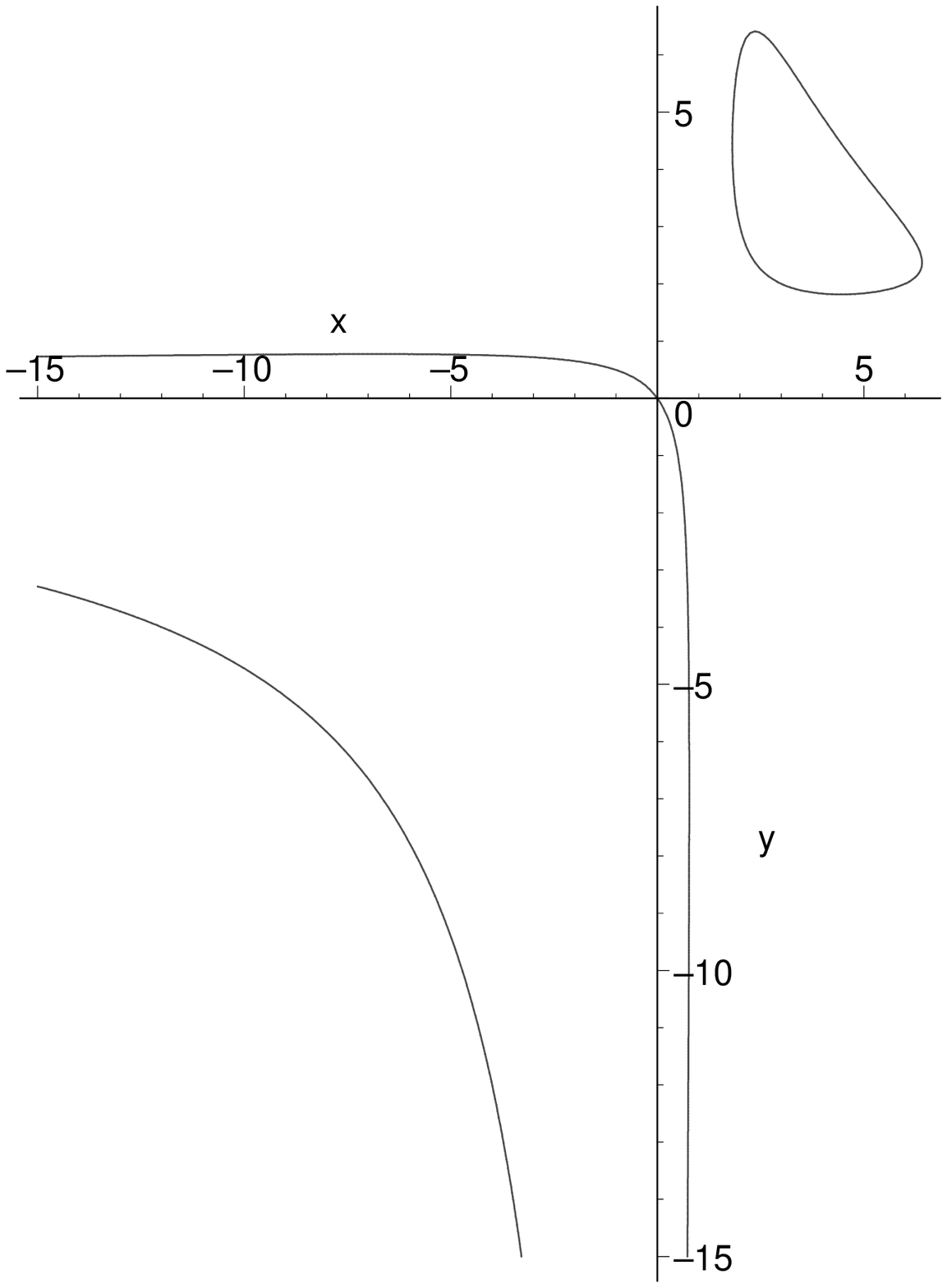}
\hskip 3cm 
\epsfxsize=63mm \epsfysize=6cm
\epsf{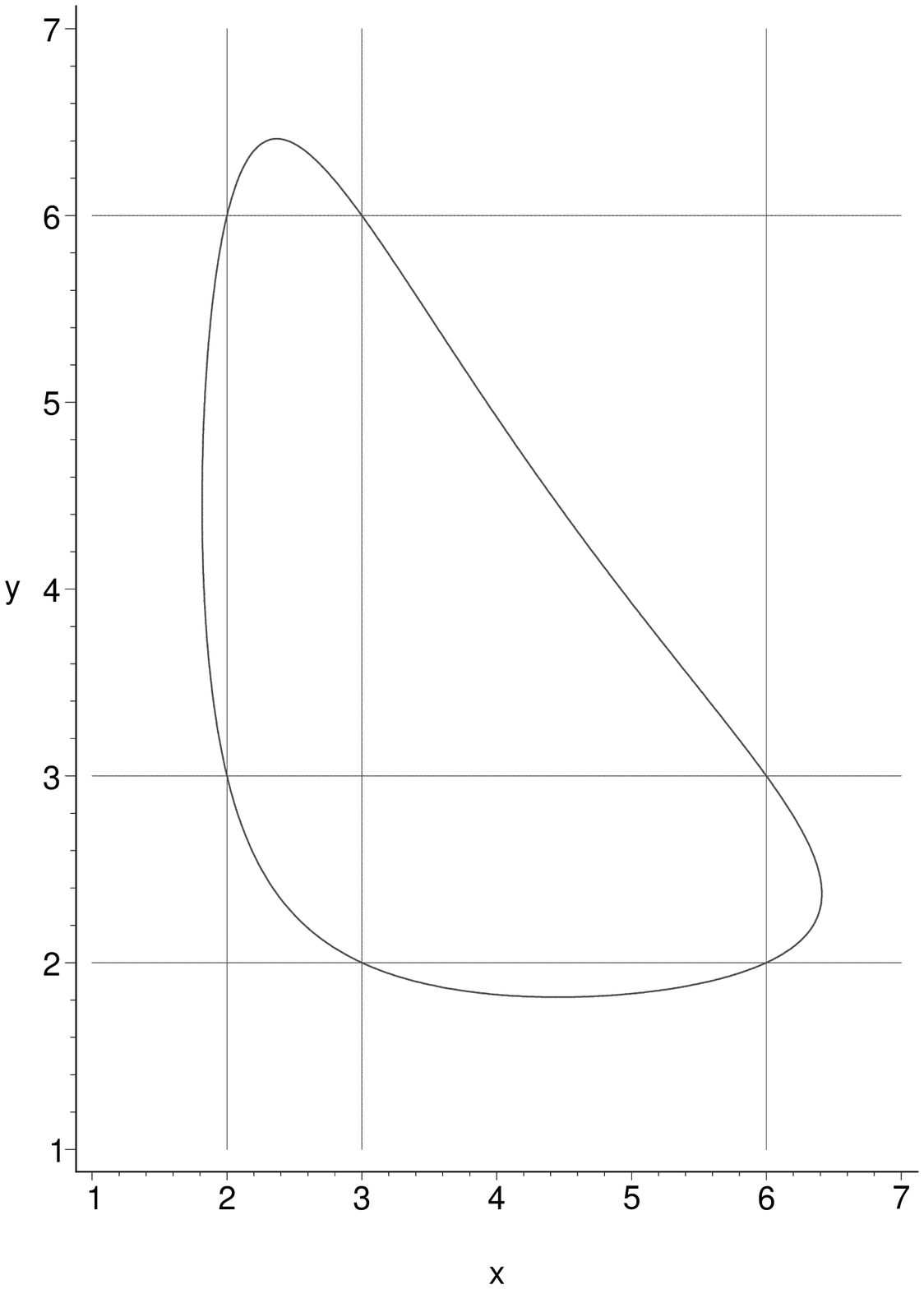}
\end{center}
\caption{The real branches of the functions $Y_i$, for $p=1/3,q=1/2$
and $r=1/6$.}
\label{figure-Y_i}
\smallskip
\noindent\hrule
\end{figure}

It will be convenient to work with an equation whose kernel is
symmetric in $x$ and $y$.  In Eq.~\Ref{eq-flatto}, let us
replace $x$ by $px$ and $y$ by $qy$. Multiplying through by $xy$ gives
\beq
(xy-x-y-pqrx^2y^2)Q(px,qy)+ x(1-qy) Q(px,0) + y(1-px) Q(0,qy)=0.
\label{eq-flatto-sym}
\eeq
Let $K(x,y)=xy-x-y-pqrx^2y^2$ denote
the kernel of this equation. The discriminant of $K$, taken as a polynomial in $y$, is
$(x-1 )^2-4pqrx^3$. Let us denote
$$
\Delta(x)= (1-\bx )^2-4pqrx.
$$
\begin{Lemma}
\label{lemma-Delta}
  Assume  $r<\min(p,q)$. Let $m=\min(p,q)$ and
$M=\max(p,q)$.
The three roots of  $\Delta(x)$, denoted $x_i$,
$i=0,1,2$,  are real, and satisfy
$$
0<x_0 <1 <x_1 < \frac 1 M \le \frac 1 m < \frac 1
r  \le x_2.
$$
\end{Lemma}
{\bf Proof.}
The variations of $\Delta(x)$ are
easy to study. Note that $\Delta(x)$ is a square when $x=1/p,1/q$, or
$1/r$.
\cqfd
Recall that $w$ is defined as  the smallest positive solution of
$w=2+pqrw^3.$ This implies that $1/(pqrw^2)=x_2$. The lemma above
gives the following factorization of $\Delta$, which will play the
role of the {\em canonical factorization\/} of Section~\ref{section-count}: 
\beq
\label{factorization}
\Delta(x)=\Delta_0 \Delta_+(x) \Delta_-(\bx)
\eeq
with
\beq
\label{factorization-values}
\Delta_0= 4pqrx_2=4/w^2, \quad \Delta_+(x)= 1 -x/x_2=1-pqrw^2x, \quad
\Delta_-(\bx)= (1-\bx x_0)(1-\bx x_1).
\eeq
As a polynomial in $y$, the kernel $K(x,y)$ of Eq.~\Ref{eq-flatto-sym}
has two roots:
$$
Y_0(x)=\displaystyle \frac{1-\bx -\sqrt{\Delta(x)}}{2pqrx},
\quad
 Y_1(x)=\displaystyle \frac{1-\bx +\sqrt{\Delta(x)}}{2pqrx}.
$$
The elementary symmetric functions of the $Y_i$ are polynomials in
$\bx =1/x$: 
$$
Y_0+Y_1= \frac {\bx (1-\bx)}{pqr}  \quad \hbox{and}
\quad Y_0Y_1 = \frac \bx {pqr} .
$$
Using the canonical factorization of $\Delta$, we see that $Y_0$ and $Y_1$
are at least analytic in the annulus $x_1<|x|<x_2$.
Let us study these functions a bit more precisely when $x$ is real.
The following lemma is illustrated by Figure~\ref{figure-Y_i}.

\begin{Lemma}
\label{lemma-Y01}
We still assume that $r<\min(p,q)$.   The functions
$Y_0$ and $Y_1$ are well-defined and real for $x \in
(-\infty, x_0] \cup [x_1, x_2]$. In particular,
$$
Y_0(1/M)=1/m, \quad Y_0(1/m)=Y_0(1/r)=1/M, \quad
Y_1(1/M)=Y_1(1/m)=1/r, \quad Y_1(1/r)=1/m.
$$
Each of the derivatives $Y'_0(x)$ and $Y_1'(x)$ admits a unique zero
on the interval $[x_1,x_2]$, respectively  denoted by $v_2$ and
$v_1$. Moreover, 
$$
x_1<\frac 1 M \le v_1=w \le \frac 1 m <v_2 < \frac 1 r \le x_2.$$ The
function $Y_0$ decreases between $x_1$ and $v_2$, and increases
between $v_2$ and $x_2$, while the function $Y_1$ increases between
$x_1$ and $v_1$, and then decreases up to $x_2$.\\
Finally, for $x \in(1/m,x_2)$, one has
$$
\frac 1 {pqx} < Y_0(x).
$$
\end{Lemma}
{\bf Proof.} The proof is a bit tedious, but elementary. We merely sketch
the different steps.

The first assertion comes from the study of the discriminant $\Delta(x)$ 
 (Lemma~\ref{lemma-Delta}). The values of $Y_i$ at the points $1/M,
 1/m$ and $1/r$ are obtained by a direct calculation. 

Let us now focus
 on the interval $[x_1,x_2]$. Given that $(1-\bx) Y_i=1+pqrxY_i^2$
 and $x\ge x_1 >1$, we have $0< Y_0(x) \le Y_1(x)$.
The derivatives of $Y_0$ and $Y_1$ with respect to $x$ can be written:
$$
Y_0'(x)= \frac 1 {x\sqrt{\Delta}} \left( (x-2) \frac {Y_0} x -1
\right), \quad
Y_1'(x)= -\frac 1 {x\sqrt{\Delta}} \left( (x-2) \frac {Y_1} x -1 \right).
$$
Given that $Y_0$ and $Y_1$ are positive on  $[x_1,x_2]$, any root of
these derivatives will be larger than~$2$.  The equation
satisfied by the $Y_i$ implies that these roots are also solutions of
$x-2=pqrx^3$. The polynomial $pqrz^3-z+2$ has two roots larger than
$2$. Let us denote them $v_1$ and $v_2$, with $v_1<v_2$. Note that
$v_1$ is actually the number $w$ defined in
Theorem~\ref{theorem-stationary}. If $x=v_i$, then
$\Delta(x)= (x-3)^2$. Hence $v_i$ belongs to the interval
$[x_1,x_2]$. Evaluating the numerators of $Y_0'$ and $Y_1'$
at $v_1$ and $v_2$ shows that 
$v_1$ cancels $Y_1'$, while $v_2$ cancels $Y_0'$. Finally, we compute 
 the numerators of $Y_0'$ and $Y_1'$ at $x_1$ and $x_2$. This
 determines the sign of these derivatives and completes the study of
 the variations of $Y_0$ and $Y_1$.

The last assertion is proved by studying the function $y \mapsto
K(x,y)$, for $x$ fixed.
\cqfd

Let $x \in (x_1,x_2)$ and let $y$ belong to the annulus $\{Y_0(x) <
|y|<Y_1(x)\}$. Let $K_r= K/(xy)= 1-\bx -\by -pqrxy$ be the rational
version of the kernel. Then the following expansion is convergent:
\beq
\label{Kr-expansion}
\frac 1 {K_r} = \frac 1 {\sqrt{\Delta(x)}} \left( \frac 1 {1-\by Y_0} +
\frac 1 {1-y /Y_1}-1\right)
=  \frac 1 {\sqrt{\Delta(x)}} \left( \sum_{n\ge 0} \by ^n Y_0^n +
 \sum_{n\ge 1} y ^n Y_1^{-n}\right).
\eeq
Let $F(y)$ be analytic in the same annulus. Let us write
$$
F(y)= \sum_{i\in \zss} f_i y^i = F^-(\by)+f_0+F^+(y),
$$
where $F^+$ and $F^-$ are the positive and negative parts of $F$.
Then 
\beq
\label{coeff-constant}
\begin{array}{llll}
[y^0] \displaystyle \frac{F(y)}{K_r}& =&
\displaystyle \frac 1 {\sqrt{\Delta(x)}} \left( F^-(1/Y_1)+f_0+F^+(Y_0)\right)\\
& =&
\displaystyle \frac 1 {\sqrt{\Delta(x)}} \left( F^-(pqrxY_0)+f_0+F^+(Y_0)\right).
\end{array}
\eeq
In particular, if $F(y)=G(y)-G(\bx \by /(pqr))$, then 
\beq \label{coeff-constant-nul}
[y^0] \displaystyle \frac{F(y)}{K_r} =0.
\eeq

\subsubsection{Simplification of the functional equation}
\label{section-simplification}
Let us go back to Eq.~\Ref{eq-flatto-sym}. 
The function $\bar Q(x,y) =
Q(px,qy)/((1-px)(1-qy))$ satisfies an equation that is symmetric in
$x$ and $y$:
\beq
(xy-x-y-pqrx^2y^2)\bar Q(x,y)+ x\bar Q(x,0)+y\bar Q(0,y)=0.
\label{eq-flatto-sym-Qbarre}
\eeq
Yet, we shall see that this function is,
in general, not symmetric in $x$ and $y$
(Section~\ref{section-comments}). We shall symmetrize it by 
considering $\bar Q(x,y)+\bar Q(y,x)$. More precisely, we shall study
separately the functions $S(x,y)$ and $D(x,y)$ (as in Sum and
Difference) defined by
\beq
\label{SDdef}
\begin{array}{lll}
S(x,y) &=& (1-qx)(1-py)Q(px,qy) +  (1-px)(1-qy)Q(py,qx),\\
D(x,y) &=& (1-qx)(1-py)Q(px,qy) -  (1-px)(1-qy)Q(py,qx).
\end{array}
\eeq
By Lemma~\ref{lemma-convergence}, these functions are analytic in
\beq
\label{domaine}
{\cal D} =
\left\{ |x|<\frac 1 r , |y| < \frac 1 M\right\} \cup 
\left\{ |x|< \frac 1 M , |y| < \frac 1 r\right\}
\eeq
with $M=\max(p,q)$. Note that $S(x,y)$ and $D(x,y)$ satisfy {\em the
  same equation\/}:
\beq
\label{S-D-eqs}
\begin{array}{lll}
S(x,y)(xy-x-y-pqrx^2y^2)+x(1-py)(1-qy)S(x,0)+y(1-px)(1-qx)S(0,y)=0,\\
D(x,y)(xy-x-y-pqrx^2y^2)+x(1-py)(1-qy)D(x,0)+y(1-px)(1-qx)D(0,y)=0.
\end{array}
\eeq
However, $S(x,y)$ is symmetric in $x$ and $y$, while
$D(x,y)=-D(y,x)$. We shall solve separately the two equations, taking
into account the respective symmetry or anti-symmetry condition. For each
equation, we will obtain  a unique solution, up to a multiplicative
factor.  A relation between the two factors will  be found
by noticing that, by definition of $S$ and $D$,
\beq
\label{double-noyau}
S\left( \frac 1 q ,0\right )+D\left( \frac 1 q ,0\right )=0.
\eeq

First  we  rewrite the above equations on $S$ and
$D$, by expressing $x(1-py)(1-qy)$ in terms of the kernel: indeed,
$$
rx(1-py)(1-qy)= (r-\bx)(x+y-xy(p+q))-\bx K(x,y).
$$
Let $K_r= K/(xy)= 1-\bx -\by -pqrxy$. The equation satisfied by $S$
can be rewritten as
\beq
\label{eq-ST}
\begin{array}{lll}
\displaystyle K_r\ \frac{xyrS(x,y)-yS(x,0)-xS(0,y)}{x+y-xy(p+q)}&=&
\displaystyle (\bx-r)S(x,0)+(\by-r)S(0,y)\\
&=& \bx T(x) +\by T(y)
\end{array}
\eeq
with 
\beq
\label{Tdef}
T(x)=(1-rx)S(x,0)=(1-rx)S(0,x).
\eeq
Similarly,
%
$$
\begin{array}{lll}
\displaystyle 
K_r\ \frac{xyrD(x,y)-yD(x,0)-xD(0,y)}{x+y-xy(p+q)}&=&
(\bx-r)D(x,0)+(\by-r)D(0,y)\\
&=& E(x) -E(y),
\end{array}
$$
where 
\beq
\label{Edef}
E(x)=(\bx-r)D(x,0)=-(\bx-r)D(0,x).
\eeq
 We have taken into
account the fact that $D(x,y)$ is anti-symmetric, so that, in
particular, $D(0,0)=0$. The functions $T$ and $E$ are analytic for
$|x|<1/r$.

From the fact that the curve $x+y=xy(p+q)$ intersects the domain of
convergence of $S(x,y)$ (near the origin), we derive from~\Ref{eq-ST} the existence of
a function $U(x,y)$, analytic in $\cal D$, such that
$$
xyrS(x,y)-yS(x,0)-xS(0,y)= (x+y-xy(p+q)) U(x,y).
$$
A similar statement holds for the function
$xyrD(x,y)-yD(x,0)-xD(0,y)$. But this function vanishes as soon as
$x=0$ or $y=0$, so that we can actually write
$$
xyrD(x,y)-yD(x,0)-xD(0,y)=xy (x+y-xy(p+q))F(x,y),
$$
for a function $F$ that is analytic is  $\cal D$. Finally, we shall need the
following initial conditions
\beq
U(x,0)=-S(0,0)=-2Q(0,0), \quad T(0)=2Q(0,0).
\label{init}
\eeq
The equations~\Ref{S-D-eqs} have thus been replaced by the following
simpler equations:
\begin{eqnarray}
K_r U(x,y)&=& \bx T(x) +\by T(y),\label{Uxy}\\
K_r xy F(x,y) &=& E(x) -E(y),\label{Fxy}
\end{eqnarray}
where $K_r=1-\bx -\by -pqrxy$  and $U,T,F, E$ are analytic in the
domain $ \cal D$ defined by~\Ref{domaine}.

\subsubsection{The algebraic kernel method}
\label{section-algebraic-analytic}
We apply the {algebraic} kernel method of
Section~\ref{section-algebraic-kernel} to the
equations~(\ref{Uxy}--\ref{Fxy}). We observe that 
$K_r$ satisfies the following invariance condition:
$$
K_r(x,y)= K_r\left( \frac{\bx \by}{pqr},y\right)
= K_r\left( x,\frac{\bx \by}{pqr}\right) \equiv K_r.
$$
Recall that $m=\min(p,q)$ and $M=\max(p,q)$. Let us fix $x$ in the
interval $(1/m,1/r)$ and restrict $y$ 
to the annulus
\beq \label{ineq1}
 {Y_0(x)} < |y| < \frac 1 M.
\eeq
By Lemma~\ref{lemma-Y01}, this annulus is non-empty, and moreover,
\beq\label{ineq2}
\frac 1 {pqx} < |y| < Y_1(x).
\eeq
 The pairs $(x,y)$ and $(\bx \by /(pqr), y)$ both belong to the
domain of convergence $\cal D$, and we thus have, in addition
to~(\ref{Uxy}--\ref{Fxy}):
\begin{eqnarray}
K_r \ U\left(  \frac{\bx \by}{pqr},y\right)&=& 
pqrxy \ T\left( \frac{\bx \by}{pqr}\right) +\by T(y),\label{Uxpy}\\
K_r  \frac{\bx }{pqr} F\left( \frac{\bx \by}{pqr},y\right) &=& 
E\left( \frac{\bx \by}{pqr}\right) -E(y).\label{Fxpy}
\end{eqnarray}
A linear combination of~\Ref{Uxy} and~\Ref{Uxpy} gives
$$
  2U(x,y) - U\left(  \frac{\bx \by}{pqr},y\right) =
\frac 1{K_r} \left( 
2\bx T(x) +\by T(y)-pqrxy T\left( \frac{\bx \by}{pqr}\right)
 \right).
$$
In view
of~(\ref{ineq1}--\ref{ineq2}), the expansion of $1/K_r$ given
by~\Ref{Kr-expansion} is convergent. Recall that $x$ is fixed; we can
now use~\Ref{coeff-constant-nul} to extract from the above equation the
coefficient of $y^0$. We obtain
\beq
  2U(x,0) - U_d\left(  \frac{\bx }{pqr}\right) =
\frac {2\bx T(x)}{\sqrt{\Delta(x)}}
\label{UT}
\eeq
where $U_d(x)$ denotes the diagonal of the series $U(x,y)$:
$$
U(x,y) = \sum_{i,j \ge 0} u_{i,j} x^i y^j \Rightarrow
U_d(x)=\sum_{i \ge 0} u_{i,i} x^i.
$$
Given that $U(x,y)$ converges absolutely in the domain $\cal D$ given
by~\Ref{domaine}, the sub-series 
$U_d(x)$ is convergent for $|x|<1/(rM)$.
By analytic continuation, \Ref{UT} holds in the annulus $\{ 1/ m
< |x|<  1 /r\}$.
Recall that $U(x,0)$ is actually a constant (see~\Ref{init}). We
now use the canonical factorization of $\Delta(x)$, given
by~\Ref{factorization}. 
We multiply~\Ref{UT} by $\sqrt{\Delta_-(\bx)}$:
$$
-\sqrt{\Delta_-(\bx)} \left( 4Q(0,0) +U_d \left(  \frac{\bx
}{pqr}\right) \right)=\frac {2\bx T(x)}{\sqrt{\Delta_0 \Delta_+(x)}}.
$$
 Using Proposition~\ref{annulus}, we can
 extract the nonnegative part of this function. Given that
 $U_d(0)=U(0,0)=-2Q(0,0)$ we obtain, using~\Ref{factorization-values}:
\beq
\label{Tsol}
T(x) = 2 Q(0,0) \left( 1 -\frac x w \right)\sqrt{\Delta_+(x)}.
\eeq
In view of~(\ref{eq-ST}-\ref{Tdef}), we have completed the
determination of the Sum function $S(x,y)$.

\medskip
Let us now work with the equations~\Ref{Fxy} and~\Ref{Fxpy}. Let us
divide~\Ref{Fxy} by $K_r$ and extract the coefficient of $y^0$. We
obtain
$$
E(x)=E(Y_0).
$$
If we do the same with~\Ref{Fxpy}, we simply find $0=0$ (recall that
$F(x,y)$ is antisymmetric). However, if
we extract instead the coefficient of $y$, we obtain,
using~\Ref{coeff-constant}: 
$$
\frac \bx{pqr} F_1\left(\frac \bx{pqr}\right) = E(0)-E(Y_0),
$$
where
$$
F(x,y) = \sum_{i,j \ge 0} f_{i,j} x^i y^j \Rightarrow
F_1(x)=\sum_{i \ge 0} f_{i,i+1} x^i.
$$
By combining both equations, and extracting the nonnegative part,  one
sees that $E(x)$ is actually a constant $E(0)$.  In view
of~(\ref{Tdef}--\ref{Edef}) and~\Ref{Tsol}, one has
$$
S(x,0)=S(0,x)= 2\frac{ Q(0,0)}{1-rx}
 \left( 1 -\frac x w \right)\sqrt{\Delta_+(x)} \quad \hbox{and} \quad
D(x,0)=-D(0,x)=\frac {xE(0)}{1-rx}.
$$
The identity~\Ref{double-noyau} completes the determination of $E$:
$$
E(0)= -2Q(0,0) \left( q - \frac 1 w\right) \sqrt{\Delta_+(1/q)}.
$$
We can now express $Q(x,0)$ and $Q(0,x)$ explicitly,
using~\Ref{SDdef}. Thanks to~\Ref{double2}, this gives exactly
Theorem~\ref{theorem-stationary}. \cqfd

\subsection{Comments on the solution}
\label{section-comments}
\subsubsection{Asymptotics}
For a good understanding of the solution of
Theorem~\ref{theorem-stationary}, or, equivalently, of the stationary
distribution of Corollary~\ref{coro-stationary}, it is useful to
determine the dominant singularities of the functions
$Q(x,0)$ and $Q(0,y)$, and hence the asymptotic behaviour of the
numbers $p_{i,0}$ and $p_{0,j}$. This is why we briefly rederive below a
result already proven in~\cite{flatto-hahn}.
\begin{Proposition}
\label{proposition-asympt}
Assume $r<\min(p,q)$. The asymptotic decay of the stationary
probabilities $p_{i,0}$ 
depends on the relative values of $p$ and $q$:
\begin{itemize}
\item If $p=q$, then $Q(x,0)$ is the reciprocal of a square root, and, 
as $i$ goes to infinity,
$$
p_{i,0} \sim c \ (r/p)^i \ i^{-1/2}
$$
for some positive constant $c$.
\item If $p<q$, then $Q(x,0)$ has a simple pole at $p/r$ as its
unique dominant singularity. The decay of the numbers $p_{i,0}$ is given by
$$
p_{i,0} \sim c \ (r/p)^i .
$$
\item If $p>q$ then $Q(x,0)$ has a square root singularity at
$1/(qrw^2)$ as its unique dominant singularity, and
$$
p_{i,0} \sim c \ (qrw^2)^i i^{-3/2}.
$$
\end{itemize}
\end{Proposition}
{\bf Proof.} We use standard results that relate the singularities of
a series to the asymptotic behaviour of its coefficients (see,
e.g.,~\cite{flajolet-odlyzko}). 

When $p=q$, the result is clear in view of
Theorem~\ref{theorem-stationary}. 
Otherwise,  the three possible singularities of
$Q(x,0)$ are   $p/q , p/r$ and $1/(qrw^2)$.
The inequalities~\Ref{ineqwpqr} imply
$$
\frac p q  < \frac p r < \frac 1 {qrw^2}.
$$
Hence our first candidate for the radius of $Q(x,0)$ is
$p/q$. However,  the numerator of $Q(x,0)$  vanishes at this
point, so that there is no pole at $p/q$. Our next candidate is
$p/r$. For this value of  
$x$, the numerator of $Q(x,0)$ is
$$
\left(1-\frac 1 {rw}\right) \sqrt{1-pqw^2} -\frac{q}{r}
 \left(1-\frac 1 {qw}\right) \sqrt{1-prw^2}.
$$
According to~\Ref{double2}, the first term
in this difference 
is negative. If  $p<q$, then  the second term is positive. Hence
the difference is 
negative, and $Q(x,0)$ has indeed a simple pole at $p/r$. 

However, if $p>q$, then~\Ref{double2} shows that the numerator of $Q(x,0)$ cancels at
$x=p/r$, so that the only singularity of $Q(x,0)$ is a square root
singularity at $1/(qrw^2)$.

Note that one can  compute explicitly, in the same way,  the multiplicative
constants denoted $c$ in the proposition.\cqfd

\subsubsection{An asymmetry of the solution}
As observed above, the function $\bar Q(x,y) =
Q(px,qy)/((1-px)(1-qy))$ satisfies an equation that is symmetric in
$x$ and $y$ (Eq.~\Ref{eq-flatto-sym-Qbarre}).
Hence  we could expect $\bar Q(x,y)$ to
be a symmetric function of $x$ and $y$. This is equivalent to the
condition
$$
\bar Q(x,0)-\bar Q(0,x)= \frac{Q(px,0)}{1-px} -\frac{Q(0,qx)}{1-qx}=0.
$$
However, we derive from Theorem~\ref{theorem-stationary} and
Eq.~\Ref{double2} that
$$
\bar Q(x,0)-\bar Q(0,x)
=\frac{2Q(0,0)x(pw-1) \sqrt{1-qrw^2}}{w(1-px)(1-qx)(1-rx)}.
$$
By~\Ref{ineqwpqr}, this quantity differs from $0$, unless $pw=1$, which forces $p=q$. If $p=q$, the solution  satisfies $Q(x,0)=Q(0,x)$, so
that the symmetry property naturally holds.
Otherwise, the asymmetry of the result comes from the asymmetric
 conditions we have 
required: $\bar Q(x,y)$ must converge when $|px|<1$ and $|qy|<1$.

\subsubsection{Flatto and Hahn's expression}
 Assume $p<q$. Eq.~\Ref{double2} shows that
that the numerator of $Q(0,y)$ vanishes when $y=q/p$ and 
$y=q/r$. 
Let us denote $\delta(y)=\sqrt{1-yprw^2}$. Then the numerator of $Q(0,y)$ is a
polynomial in $\delta(y)$, of degree 3, and two of its roots are
$\delta(q/p)=\sqrt{1-qrw^2}$ and $\delta(q/r)=\sqrt{1-pqw^2}$.  The
third root is then easily 
determined to be
$-\sqrt{1-prw^2}$. Hence, up to a multiplicative constant independent
of $y$, the numerator of $Q(0,y)$ factors as
$$
  \left(\sqrt{1-yprw^2} - \sqrt{1-qrw^2}\right)
 \left(\sqrt{1-yprw^2}-  \sqrt{1-pqw^2}\right)
 \left(\sqrt{1-yprw^2} + \sqrt{1-prw^2}\right) .
$$
The denominator of $Q(0,y)$ is already factored in $y$, and also
vanishes at $y=q/p$ and $y=q/r$.  Up to a
multiplicative constant, it factors as
$$
\left(\sqrt{1-yprw^2} + \sqrt{1-qrw^2}\right)\left(\sqrt{1-yprw^2} - \sqrt{1-qrw^2}\right)
\hskip 19cm
$$
$$
\hskip 6cm
\times \left(\sqrt{1-yprw^2} + \sqrt{1-pqw^2}\right) \left(\sqrt{1-yprw^2} -\sqrt{1-pqw^2}\right).
$$
Two simplifications occur, and finally
\beq
Q(0,y)= Q(0,0) \frac{ \Psi(y)} {\Psi(0)}
\label{flattoy}
\eeq
where
$$
\Psi(y) = \frac{ \sqrt{1-yprw^2} + \sqrt{1-prw^2} }{(\sqrt{1-yprw^2} +
\sqrt{1-qrw^2})(\sqrt{1-yprw^2} + \sqrt{1-pqw^2})}.
$$
Now, using~\Ref{double2}, the function $Q(x,0)$ can
be rewritten as
$$
Q(x,0)=\frac{Q(0,0)}{(1-qx/p)(1-rx/p)}
\left( \left(1-\frac x {pw}\right) \sqrt{1-xqrw^2} +x
 \left(1-\frac 1 {pw}\right) \sqrt{1-qrw^2} \right) .
$$
In this form, the numerator of $Q(x,0)$ now looks more like the
numerator of $Q(0,y)$. More precisely, denoting the latter numerator
by $P(\delta(y))$, the former numerator is exactly
$-P(-\sqrt{1-xqrw^2})$, and hence factors as
$$
\left( \sqrt{1-xqrw^2} + \sqrt{1-qrw^2} \right) \left(\sqrt{1-xqrw^2} +\sqrt{1-pqw^2}\right) 
\left(\sqrt{1-xqrw^2} - \sqrt{1-prw^2}\right).
$$
The denominator of $Q(x,0)$ is also easily factored; two
simplifications occur again, and we end up with
\beq
Q(x,0)= Q(0,0)\frac{ \Phi(x)}{ \Phi(0) }
\label{flattox}
\eeq
where
$$
\Phi(x) = \frac{ \sqrt{1-xqrw^2} + \sqrt{1-qrw^2} }{(\sqrt{1-xqrw^2} +
\sqrt{1-prw^2})(\sqrt{1-xqrw^2} - \sqrt{1-pqw^2})}.
$$
The expressions~\Ref{flattoy} and~\Ref{flattox} are, with our
notation, the forms given in Flatto and Hahn's
paper~\cite{flatto-hahn}. They are nicely factored, and it is easy to
derive from them the singularities of $Q(x,0)$ and $Q(0,y)$. However,
they have two drawbacks: first, they are only valid when $p\le q$, and
hide the symmetry of the result in $p$ and $q$, which is clear from the
expressions of Theorem~\ref{theorem-stationary}. Secondly, they somehow 
contain ``two many'' radicals, and suggest that $Q(x,0)$ and $Q(0,y)$
will be algebraic of degree $3\times 2^4$ over the field
$\qs(p,q,x,y)$, whereas as suggested by
Theorem~\ref{theorem-stationary}, they have only degree $3\times
2^2=12$. This can be checked using a computer algebra package, like
{\sc Maple}.

\section{Enumeration \& probability: the law of the chain}\label{section-law}
In this section, we consider again the Markov chain illustrated in
Figure~\ref{figure-transitions}. We start this chain at time $0$ at
the origin of the lattice, and 
address the question of computing the probability $p_{i,j}(n)$ that
the walk reaches the point $(i,j)$ at time $n$. This question is, in
essence, close to Section 2: we are again {\em
enumerating\/} paths
according to a certain weight. This weight is the probability that the
trajectory begins with this path. But this question is also related to
Section 3, since we 
expect the probability $p_{i,j}(3n-i-j)$ to converge to $3p_{i,j}$
as $n$ goes to infinity, when $r< \min(p,q)$, where $p_{i,j}$ is the stationary
distribution of the chain (the factor $3$ accounts for the
periodicity of the chain).

The notation we adopt is similar to that of Section 3:
we introduce the following four  \gfs \ for the probabilities $p_{i,j}(n)$:
$$
P_{0,0}=\sum_{n \ge 0} p_{0,0}(n) t^n,
 \quad P_1(x)=\sum_{n, i > 0} p_{i,0}(n) x^it^n,
\quad  P_2(y)=\sum_{n, j > 0} p_{0,j}(n) y^jt^n,
$$
$$P(x,y) =\sum_{n,i,j > 0} p_{i,j}(n) x^i y^jt^n.
$$
The step by step construction of the walks gives  the following functional
equation:
$$(1-p\bx t -q\by t -rxyt) P(x,y) +(1-p'\bx t -r' xyt) P_1(x) 
+(1-q''\by t -r''xyt) P_2(y)+(1-xyt) P_{0,0}=1.
$$
Again, we assume that the transition probabilities on the border of
the quadrant are related to those inside the quadrant by the
conditions~\Ref{conditions}. This allows us to rewrite the above
functional as
$$
( 1-p\bx t -q\by t -rxyt) Q(x,y)+ q(\by t-1) Q(x,0) + p(\bx t-1) Q(0,y) = r
$$
with
\beq
\label{Q-P}
Q(x,y) = P_{0,0}+ r'P_1(x)+r''P_2(y)+rP (x,y).
\eeq
It will be convenient to have a kernel symmetric in $x$ and $y$, and 
our starting point will actually be
\beq
(xy-t(x+y+pqrx^2y^2))Q(px,qy) + (t-qy)xQ(px,0) + (t-px)yQ(0,qy) =rxy.
\label{eqsymt}
\eeq
%
We are back to the (safe) world of formal power series in $t$ with
coefficients in $\qs(x,y)$, and we will mimic the obstinate kernel
method of Sections~\ref{section-obstinate}
and~\ref{section-symmetric}. The only new difficulty arises from the
absence of symmetry, since $Q(x,0)\not = Q(0,x)$ when $p\not = q$. \\
The  kernel of the above equation, considered
as a polynomial in $y$, has two roots,
$$
\begin{array}{lclllll}
Y_0(x)&=&\displaystyle \frac{1-t\bx -\sqrt{(1-t\bx )^2-4pqrt^2x}}{2pqrtx}&=&
 & &t+\bx t^2 + O(t^3) , \\
\\
 Y_1(x)&=&\displaystyle \frac{1-t\bx +\sqrt{(1-t\bx )^2-4pqrt^2x}}{2pqrtx}
&=&\displaystyle\frac \bx {pqrt} -\frac{\bx ^2}{pqr} &-&t-\bx t^2 + O(t^3) .
\end{array}$$
The elementary symmetric functions of the $Y_i$ are again polynomials
in $1/x$:
$$
Y_0+Y_1= \frac {\bx(1-t\bx) }{pqrt} \quad \hbox{and}
\quad Y_0Y_1 = \frac{\bx}{pqr} .
$$
The discriminant $\Delta(x)=(1-t\bx )^2-4pqrt^2x$ vanishes for three
values of $x$: two of them, say 
$X_0$ and $X_1$,  are power series in $\sqrt t$ while the third one,
$X_2$, is a Laurent series in $t$ that starts with a term in
$t^{-2}$. Let us define $Z\equiv Z(t)$ to be  the unique power series
in $t$ such that 
$$Z=1+4pqrt^3Z^3.$$
 Then $4pqrt^2X_2Z^2=1$, and  the
canonical factorization of $\Delta(x)$ reads 
$$
\Delta(x)=\Delta_0 \Delta_+(x) \Delta_-(\bx)
$$
with
\beq
\Delta_0= 4pqrt^2 X_2 = \frac 1 {Z^2},
 \quad \Delta_+(x)= 1 -x/X_2 = 1-4pqrt^2Z^2x ,
\label{Delta0+}
\eeq
\beq
 \Delta_-(\bx)= (1-\bx X_0)(1-\bx X_1) = 1-tZ(1+Z)\bx + t^2Z^2 \bx^2. 
\label{Delta-}
\eeq
As in Section~\ref{section-stationary},	it will be convenient to
handle two functions $S(x,y)$ and $D(x,y)$,  which are respectively
symmetric and antisymmetric in $x$ and $y$. We define them by
\beq
\label{SDdef-t}
\begin{array}{lll}
S(x,y) &=& (t-qx)(t-py)Q(px,qy) +  (t-px)(t-qy)Q(py,qx),\\
D(x,y) &=& (t-qx)(t-py)Q(px,qy) -  (t-px)(t-qy)Q(py,qx).
\end{array}
\eeq
Then
$$t\left(xy-t(x+y+pqrx^2y^2)\right)S(x,y)+
\hskip 14cm 
$$
\beq\label{S-eq-t}
\hskip 3cm (t-py)(t-qy)xS(x,0)+(t-px)(t-qx)yS(0,y)=\F(x,y)+\F(y,x),
\eeq
$$
t\left(xy-t(x+y+pqrx^2y^2)\right)D(x,y)+
\hskip 14cm 
$$
\beq 
\label{D-eq-t}
\hskip 3cm (t-py)(t-qy)xD(x,0)+(t-px)(t-qx)yD(0,y)=\F(x,y)-\F(y,x),
\eeq
where
$$
\F(x,y) = rxyt(t-qx)(t-py).
$$
\subsection{Statement of the results}
After all the
algebraic series we have met, one might expect the probability
\gf \ of the law of the chain to be algebraic again. This is, however,
only true if $p=q$. 
\begin{Theorem}[The symmetric case]
\label{thm:symmetric}
Assume $p=q$. The  three-variate \gf \ for the probabilities
$p_{i,j}(n)$ is algebraic, and can be expressed explicitly in terms
of the unique power series $Z\equiv Z(t)$ satisfying
$Z=1+4pqrt^3Z^3$.  In particular, the \gf \ of walks ending at
the origin is algebraic of degree $6$:
$$
P_{0,0}=\sum_{n \ge 0} p_{0,0}(3n) t^{3n} = 
\frac r p \left ( \frac {\sqrt{1-pZ(1+Z)+p^2Z^2}}{1-2pZ}-1\right)
= \frac r p \left (\frac {\sqrt{\Delta_-(p/t)}}{1-2pZ}-1\right) ,
$$
where $\Delta_-(\bx)$ is given by~{\em\Ref{Delta-}}.
More generally, the series $Q(px,0)= P_{0,0}+ r' P_1(px)$ is given by:
$$
(t-x(1-p)+prx^2t^2)Q(px,0)= \frac r {2p} \left(
\frac  { (2tZ-x)\sqrt{\Delta_-(p/t)}\sqrt{\Delta_+(x)}}{Z(1-2pZ)}
-2t+x(1-p) \right)
$$ 
where $\Delta_+(x)$ is given by~{\em\Ref{Delta0+}}. The expression of
$P_{0,0}$ can be recovered from the value of $Q(px,0)$ by setting $x=0$.
\end{Theorem}
This theorem, and {all the results of
    this section, will be proved in Section~\ref{section-proofs}. 

What happens in the general case? We have expressed the series
 $Q(px,qy)$ in terms of two series $S(x,y)$ and $D(x,y)$, which are
 respectively symmetric and antisymmetric in $x$ and $y$.  It turns
 out that 
the Sum series $S(x,y)$ is always algebraic, while the Difference series
 $D(x,y)$ is transcendental (unless $p=q$). The algebraicity of
 $S(x,y)$  has an interesting consequence: The \gf \
 $P_{0,0}$ that counts walks ending at the origin is {\em always algebraic\/},
 even when $p\not = q$. 
\begin{Theorem}[The general case: algebraic part]\label{thm:alg-part}
The series $S(x,y)$ defined by~{\em \Ref{SDdef-t}} is algebraic and can be
expressed explicitly in terms
of the unique power series $Z\equiv Z(t)$ satisfying
$Z=1+4pqrt^3Z^3$. 
In particular, the coefficient of $x^0y^0$ in $S(x,y)$ is an algebraic series
in $t$. It is equal to $2t^2 P_{0,0}$, where $P_{0,0}$ counts walks
ending at the origin, and we have:
$$
2pqP_{0,0}+r(1-r)= A_{p,q}+A_{q,p}
$$
where $A_{p,q}$ is the following algebraic series in $t$:
$$
A_{p,q}= 
\frac{(p(1-2p)-qrt^3) \sqrt{\Delta_-(p/t)}}{(1-t^3)(1-2pZ)  }
$$
and $\Delta_-(\bx)$ is given by~{\em \Ref{Delta-}}. 
 The algebraic series
$P_{0,0}$ has degree $6$ if $p=q$, and degree $12$ otherwise.
More generally,
the series $S(x,y)$ satisfies:
$$
\left(t-(1-p)x+t^2qrx^2 \right)\left(t-(1-q)x+t^2pr x^2 \right)
\frac{S(x,0)}t+ \frac{rH(x)}{2pq}=
\hskip 10cm
$$
$$
\hskip 8cm
\frac{(2tZ-x)\sqrt{\Delta_+(x)}}{2pqZ}
\Big( A_{p,q}\G_{p,q}(x) + A_{q,p}\G_{q,p}(x)\Big),
$$
where 
$\Delta_+(x)$ is given by~{\em\Ref{Delta0+}}, and
 $\G_{p,q}(x)$ and $H(x)$ denote the following polynomials in $t$ and
$x$:
$$
\G_{p,q}(x)= (t-xq) \left(t-(1-q)x+t^2pr x^2 \right),
$$
\beq \label{Hdef}
H(x)=q(2t-x+px)\G_{p,q}(x)+p(2t-x+qx)\G_{q,p}(x)- (p-q)^2x^2
\Big( 2rt-(1-p)(1-q)x+2pqrx^2t^2\Big).
\eeq
 The expression of
$P_{0,0}$ can be recovered from the value of $S(x,0)$ by setting
 $x=0$. An expression for $S(x,y)$ can be obtained using~{\em\Ref{S-eq-t}}.
\end{Theorem}
 This theorem will allow us to complete the proof of the following
 result, announced in Section~\ref{section-stationary}.
\begin{Corollary} \label{coro-ergodic}
The Markov chain schematized in
Figure~\ref{figure-transitions}, with the border
conditions~{\em\Ref{conditions}}, is ergodic (that is, has a stationary
distribution) if and only if $r<\min(p,q)$.
\end{Corollary}
Theorem~\ref{thm:alg-part}
specializes to Theorem~\ref{thm:symmetric} when $p=q$. 
It states that walks ending at the origin have an algebraic \gf .  
What about   the \gfs \ $P_1(x)$ and
$P_2(y)$ that count  walks ending on the $x$- or $y$-axis? By
symmetry of the model, $P_1(x)$ is algebraic if and only if
$P_2(y)$ is  
algebraic too. In view of Eqs.~\Ref{Q-P},~\Ref{eqsymt} and~\Ref{SDdef-t},
this holds if and only if 
$S(x,0)$ and $D(x,0)$ are algebraic. If $p=q$, then $D(x,y)$ is obviously
zero, and Theorem~\ref{thm:symmetric} tells us
that all the \gfs\
under consideration are algebraic. If $p\not = q$, we shall prove that
$D(x,0)$ is transcendental (but D-finite), and give an explicit
expresion of it.   

\medskip
So far, we have expressed many of our series in terms of the canonical
factorization of the discriminant ${\Delta(x)}$. This is the
case, for instance, in Theorem~\ref{thm:alg-part} above, where the
expression of $S(x,0)$ involves $\sqrt{\Delta_+(x)}$, which we could
call the positive \emph{multiplicative} part of $\sqrt{\Delta(x)}$. In
order to express $D(x,0)$, we need to introduce the positive \emph{additive}
part of $\sqrt{\Delta(x)}$, as defined by~\Ref{positive-part}. More
precisely, the expression of $D(x,0)$ will involve the positive (additive) part of
\beq
B(x):= (Y_0-Y_1) (2t-x+pqrtx^3) = \frac{\sqrt{\Delta(x)}}{pqrt}
\left(1-2t\bx -pqrtx^2\right).
\label{Tdef-new}
\eeq
We shall compute below the
expansion of $B$ in $t$ and $x$, using the Lagrange inversion formula.
In particular, we will see that the positive (additive) part of $B(x)$ reads
\beq
B^+(x)= -xt -x^2 + 2 C^+(x)
 \label{TUplus}
\eeq
where all terms in $C^+(x)$ are multiples of $x^3$.
We then define the series $C^-(\bx)$ by 
\beq
B(x)= 
\frac{1-2\bx t}{pqt} + 2C^-(\bx)  -xt -x^2 + 2 C^+(x) .
\label{U-def}
\eeq
Observe that $\Delta(t/p)$, and  hence $B(t/p) $, is a well-defined
Laurent series in $t$. Clearly, $C^+(t/p)$ is well-defined too:
consequently, by difference, we can define  $C^-(p/t)$ as a Laurent
series in $t$, even though it
would be meaningless to replace $x$ by $p/t$ in the expansion of  $C^-(x)$. 
\begin{Theorem}[The general case: transcendental part]
\label{thm:trans-part}
 When $p\not = q$, the
series $D(x,y)$ defined by~{\em\Ref{SDdef-t}} is D-finite but  transcendental.
The same holds for its specialization $D(x,0)$. Consequently, the
series  $P_1(x)$ and  $P_2(y)$ which count walks 
ending on the $x$- or $y$-axis are transcendental.
The series $D(x,0)$ satisfies:
$$
\left(t-(1-p)x+t^2qrx^2 \right)\left(t-(1-q)x+t^2pr x^2 \right)
\frac{D(x,0)}{rxt} 
+x(p-q)\left(t^2(1-r)rx^2-x/2+t\right) =
\hskip 10cm
$$
$$
\hskip 5cm
rx(p-q)t^2 C^+(x)
- \frac{t}{1-t^3}
\Big( pC^-(p/t)\G_{p,q}(x)- qC^-(q/t)\G_{q,p}(x)\Big),
$$
where $C^+(x)$ and $C^-(\bx)$ are defined by~{\em
  (\ref{TUplus}--\ref{U-def})} and, as in Theorem~{\em \ref{thm:alg-part}},
$$
\G_{p,q}(x)= (t-xq) \left(t-(1-q)x+t^2pr x^2 \right).
$$
An expression of $D(x,y)$ can then be obtained using~{\em\Ref{D-eq-t}}.
\end{Theorem}
Note the similarities between the expressions of $S(x,0)$
(Theorem~\ref{thm:alg-part}) and
$D(x,0)$ (Theorem~\ref{thm:trans-part}). One could take the sum and difference of these expressions to
recover the series $Q(px,0)$ and $Q(0,qx)$, but, as no significant
simplification arises, we shall not do this. 

There is still one natural question that is not answered by the
combination of the above two theorems: we have seen that the \gf \
$P_{0,0}$ of walks ending at the origin is algebraic, but that
the series $P_1(x)$ that counts walks ending on the $x$-axis is
transcendental. Yet, for $i>0$, the 
coefficient of $x^i$ in $P_1(x)$, being
$$
P_{i,0}:=\sum_{n\ge 0} p_{i,0}(n) t^n,
$$
counts walks ending at $(i,0)$  and might be algebraic. The following
corollary tells us that this is not (systematically) the case.
\begin{Corollary}\label{corollary-Pi0}
 Some of the series $P_{i,0}$ are transcendental.
\end{Corollary}

\noindent{\bf Note.} We can obtain an explicit expression of the series
$C^+(x)$ 
and $C^-(\bx)$ 
by expanding $B(x)$  in $x$ and
$t$. Let us write
$$
%
Y_0-Y_1= 2Y_0 - (Y_0+Y_1) =2Y_0-\frac{\bx(1-\bx t)}{pqrt},
$$
and observe that the series $Y_0$, which cancels the kernel, is Lagrangian
in $t$:
$$
Y_0=t\left(1+\bx Y_0+pqrxY_0^2\right).
$$
The Lagrange inversion formula  yields
$$B(x)= 
\frac{(1-\bx t)(1-2\bx t)}{pqrt} -xt -x^2 
+2\sum_{n\ge 2}t^n \sum_{k=0}^{\lfloor{n/2}\rfloor} x^{3k-n+2}
(pqr)^k \frac{(3k-n+1)(3k-n)(n-2)!}{k!(k+1)!(n-2k)!}
.
$$
Consequently,
$$
C^+(x)=\sum_{n\ge 2}t^n \sum_{k=\lceil (n-1)/3\rceil}^{\lfloor
n/2\rfloor} x^{3k-n+2} 
(pqr)^k \frac{(3k-n+1)(3k-n)(n-2)!}{k!(k+1)!(n-2k)!}.
$$
One may also write an explicit expansion of $C^-(\bx)$.
\subsection{Proofs}\label{section-proofs}
\noindent{\bf Proof of Theorem~\ref{thm:alg-part}.} Let us start from
the equation~(\ref{S-eq-t}) defining $S(x,y)$. 
As in Section~\ref{section-obstinate}, the pairs $(x,Y_0)$ and
$(Y_0,Y_1)$ cancel the kernel 
and can be substituted for $(x,y)$ in this equation.
We thus obtain two equations:
\beq
\left\{
\begin{array}{lll}
(t-pY_0)(t-qY_0)T(x)+(t-px)(t-qx)T(Y_0)&=&\F(x,Y_0)+\F(Y_0,x), \\
(t-pY_1)(t-qY_1)T(Y_0)+(t-pY_0)(t-qY_0)T(Y_1)&=&\F(Y_0,Y_1)+\F(Y_1,Y_0), \\
\end{array}
\right.
\label{eq1}
\eeq
with $T(x)=xS(x,0)$.
Let us form a  symmetric function of $Y_0$ and $Y_1$ based on a
divided difference: We multiply the first equation by $2(t-pY_1)(t-qY_1)$
and the second one by $ (t-px)(t-qx)$, and  take the difference
of the resulting equations:
$$
2(t-pY_0)(t-qY_0)(t-pY_1)(t-qY_1)T(x)+(t-px)(t-qx)
\Big((t-pY_1)(t-qY_1)T(Y_0)-(t-pY_0)(t-qY_0)T(Y_1)\Big)
$$
$$
=2(t-pY_1)(t-qY_1) \Big(\F(x,Y_0)+\F(Y_0,x)\Big) - (t-px)(t-qx)\Big(\F(Y_0,Y_1)+\F(Y_1,Y_0)\Big).
$$
Recall that
$$
K(x,y)= xy-t(x+y+pqrx^2y^2)= -pqrtx^2(y-Y_0)(y-Y_1).
$$
We use this identity  to express the
coefficient of $T(x)$ as a Laurent polynomial in $x$ and $t$. Then, we
separate the symmetric and anti-symmetric 
parts of the right-hand side using
\begin{eqnarray*}
\Psi(Y_0,Y_1)&=& \frac 1 2 \Big( \Psi(Y_0,Y_1)+\Psi(Y_1,Y_0)\Big) 
+ \frac 1 2 \Big( \Psi(Y_0,Y_1)-\Psi(Y_1,Y_0)\Big). 
\end{eqnarray*}
This gives
$$
2\frac{(t-(1-p)x+t^2qrx^2)(t-(1-q)x+t^2pr x^2)}{pqr^2x^4}T(x) + 
\frac{t\ H(x)}{p^2q^2rx^3} = \hskip 10cm
$$
\beq
-(t-px)(t-qx)
\Big((t-pY_1)(t-qY_1)T(Y_0)-(t-pY_0)(t-qY_0)T(Y_1)\Big)
+(Y_0-Y_1) \frac{t^2 J(x)}{pqx}
\label{P-main}
\eeq
where $H(x)$ is given by~\Ref{Hdef} and $J(x)$ is also a polynomial in $x$ and $t$:
$$
J(x)=q\G_{p,q}(x)+p\G_{q,p}(x)+x(p-q)^2(t-rx).
$$
As we are getting used to the method, let us merge the next two steps:
instead of first dividing by $(Y_0-Y_1)$ and then multiplying by
$\sqrt{\Delta_-(\bx)}$, let us divide~\Ref{P-main} by
$2\sqrt{\Delta_+(x)}/(pqrx)= 
-2t(Y_0-Y_1)/\sqrt{\Delta_0 \Delta_-(\bx)}$. We obtain, in view
of~\Ref{Delta0+}: 
$$
\frac{(t-(1-p)x+t^2qrx^2)(t-(1-q)x+t^2pr x^2)}{rx^3\sqrt{\Delta_+(x)}}
T(x) + \frac{t\ H(x)}{2pqx^2\sqrt{\Delta_+(x)}} = \hskip 50mm
$$
$$\hskip 10mm
\frac{\sqrt{\Delta_-(\bx)}}{2Z}\left((t-px)(t-qx)
\frac{(t-pY_1)(t-qY_1)T(Y_0)-(t-pY_0)(t-qY_0)T(Y_1)}{t(Y_0-Y_1)}
- \frac{t\ J(x)}{pqx}\right).
$$
The left-hand side of this equation, as a Laurent series in $x$, has
valuation $-2$, while the right-hand side only involves powers of $x$
smaller than or equal to $2$. Extracting the positive part in $x$, and
mutiplying by $x^2$ gives
\beq
\frac{(t-(1-p)x+t^2qrx^2)(t-(1-q)x+t^2pr x^2)}{rx\sqrt{\Delta_+(x)}}
T(x) + \frac{t\ H(x)}{2pq\sqrt{\Delta_+(x)}}  = L(x),
\label{eq-resultat}
\eeq
where $L(x)$ is a polynomial in $x$, of degree $4$,  with coefficients in
$\qs[t,T_1,T_2,T_3]$
where $T_i\equiv T_i(t)$ denotes the coefficient of $x^i$ in $T(x)$. 
We do not give the explicit expression of $L(x)$, but refer the reader
to his/her favourite computer algebra system.

We now have to determine the {three}
unknown functions  $T_1,T_2$ and $T_3$. Fortunately, we can
compute $T(x)$ at three values of $x$ using~\Ref{eq1}: first at
$x=t/p$, then at $x=t/q$, and finally at $x=W$, where $W$ is the
unique power series in $t$ that satisfies  $K(W,W)=0$ (so that
$W=Y_0(W)$). Remarkably,  $W$ is simply
related to the parameter $Z$ by $W=2tZ$. The three values of $T(x)$
that we obtain are:
\begin{eqnarray}
T(t/p)&=&\frac{t^3(p-q)\left(q-r+\sqrt{(1-p)^2-4t^3qr}\right)}{2p^2q(1-t^3)},
\label{Ptp}\\
T(t/q)&=&\frac{t^3(q-p)\left(p-r+\sqrt{(1-q)^2-4t^3pr}\right)}{2pq^2(1-t^3)},
\label{Ptq}\\
T(W)&=&rtW^2.\nonumber
\label{eqPR}
\end{eqnarray}
Setting $x=W$ in~\Ref{eq-resultat}, we find that the left-hand side
vanishes. Hence $L(W)=0$, and this gives an expression of
$T_2$ in terms of $T_1$:
$$
T_2=-rt+\frac{2r-1}{2t}T_1.
$$
The polynomial $L(x)$ now takes the following form:
\beq
L(x)= \frac{(W-x)(t(t-x-tpqrx^3)M_0+M_1    x^2)}{pqrW}
\label{D-R}
\eeq
where 
$$
M_0= pqT_1+r(p+q)t^2,
$$
and $M_1$ involves both $T_1$ and $T_3$. 
It remains to evaluate Eq.~\Ref{eq-resultat} at $x=t/p$ and
$x=t/q$, using the expressions of $T(t/p)$ and $T(t/q)$ given
by~\Ref{Ptp} and~\Ref{Ptq}, to obtain 
$$
M_0=\frac{t^3}{1-t^3} \left( 
\frac{(p(1-2p)-t^3qr)\sqrt{\Delta_-(p/t)}}{t-pW}
+
\frac{(q(1-2q)-t^3pr)\sqrt{\Delta_-(q/t)}}{t-qW}
\right).
$$
$$
M_1=\frac{t^3}{1-t^3} \left( 
\frac{(p(1-2p)-t^3qr)(q(1-q)+prt^3)\sqrt{\Delta_-(p/t)}}{t-pW}\right.
\hskip 3cm
$$
$$
\hskip 8 cm 
+\left. 
\frac{(q(1-2q)-t^3pr)(p(1-p)+qrt^3)\sqrt{\Delta_-(q/t)}}{t-qW}
\right).
$$
Theorem~\ref{thm:alg-part} follows, using~\Ref{D-R}
and~\Ref{eq-resultat}. 
\cqfd

\bigskip

\noindent
{\bf Proof of Corollary~\ref{coro-ergodic}.} We have already seen that
the condition $r <\min(p,q)$ is necessary for the chain to have a
stationary distribution (Lemma~\ref{lemma-convergence}). Assume this
condition holds. As the chain is irreducible, it suffices to prove
that the point $(0,0)$ is positive recurrent, that is, that the
probability $p_{0,0}(3n)$ converges to a positive constant as $n$ goes
to infinity~\cite{chung}. The generating function of these numbers,
denoted $P_{0,0}$, is given explicitly in Theorem~\ref{thm:alg-part}.  

This leads us  to determine the smallest singularity of the series
$A_{p,q}$. The technique is standard for algebraic functions, and we
only sketch the main steps. 
The series $Z$ becomes singular at
$t^3=1/(27pqr)>1$.  Then, we note that
$$
p(1-2p)-qrt^3= \frac{(1-2pZ)(1+(2p-1)Z(1+2pZ))}{4pZ^3},
$$
so that there is actually no pole in $A_{p,q}$ if $Z$ reaches $1/(2p)$. Moreover,
if $\Delta_-(p/t)=0$, then $\Delta(t/p)=0$. But, as $t$
increases from $0$ to $1$, $\Delta(t/p)=(1-p)^2-4qrt^3$ decreases from
$(1-p)^2$ to $(q-r)^2$ and thus does not vanish. Hence $A_{p,q}$
(and $A_{q,p}$) has its smallest singularity at $t^3=1$, and this
singularity is a simple pole. Consequently, the coefficient of
$t^{3n}$ in $A_{p,q}$ tends to a constant as
$n\rightarrow \infty$. The same holds for   $A_{q,p}$.

It remains to show that the sum of these two constants is not zero. We
are actually going to compute them explicitly: this will not only
conclude the proof of the corollary, but also  allow us to recover the value of
$p_{0,0}$ given in Corollary~\ref{coro-stationary}. First, we note that
$Z(1)=1+4pqrZ(1)^3$, and conclude that $Z(1)=w/2$, where $w$ is the
real number defined in Theorem~\ref{theorem-stationary}. Then, the definition of the
canonical factorization gives, when $t=1$, 
$$
\Delta (1/p)=(q-r)^2=\Delta_0 \Delta_+(1/p) \Delta_-(p),
$$
and so by~\Ref{Delta0+}, 
$$
\sqrt{\Delta_-(p)} = \frac{(q-r)w}{2\sqrt{1-qrw^2}}.
$$
Thus, as $t \rightarrow 1$,
$$
2pqP_{0,0} \sim \frac w{2(1-t^3)} \left(
\frac{(p(1-2p)-qr)(q-r)}{(1-pw) \sqrt{1-qrw^2}} + 
\frac{(q(1-2q)-pr)(p-r)}{(1-qw) \sqrt{1-prw^2}}
 \right).
$$
Note that $p(1-2p)-qr= (q-p)(p-r)$. Using~\Ref{double2}, we rewrite
the above identity as
$$
2pqP_{0,0} \sim \frac 1{1-t^3} \frac{w(p-r)(q-r) |p-q| }{(1-rw)
  \sqrt{1-pqw^2}}.
$$
It follows that, as $n\rightarrow \infty$,
$$
p_{0,0}(3n) \rightarrow \frac{w(p-r)(q-r) |p-q| }{2pq(1-rw)
  \sqrt{1-pqw^2}}.
$$
Given that the chain has period $3$, this agrees with
Corollary~\ref{coro-stationary}. 
\cqfd

\bigskip

\noindent{\bf Proof of Theorem~\ref{thm:trans-part}.}
 Let us start from
the equation~(\ref{D-eq-t}) defining $D(x,y)$. 
The pairs $(x,Y_0)$ and $(Y_0,Y_1)$ cancel the kernel
and can be substituted for $(x,y)$ in this equation.
We thus obtain:
$$
\left\{
\begin{array}{lll}
(t-pY_0)(t-qY_0)E(x)-(t-px)(t-qx)E(Y_0)&=&\F(x,Y_0)-\F(Y_0,x), \\
(t-pY_1)(t-qY_1)E(Y_0)-(t-pY_0)(t-qY_0)E(Y_1)&=&\F(Y_0,Y_1)-\F(Y_1,Y_0), \\
\end{array}
\right.
$$
with $E(x)=xD(x,0)=-xD(0,x)$.
We now want to form a symmetric function of $Y_0$ and $Y_1$ based on a
sum. We multiply the first equation by $2(t-pY_1)(t-qY_1)$
and the second one by $ (t-px)(t-qx)$, and take the
sum of the resulting equations: 
$$
2(t-pY_0)(t-qY_0)(t-pY_1)(t-qY_1)E(x)-(t-px)(t-qx)
\Big((t-pY_1)(t-qY_1)E(Y_0)+(t-pY_0)(t-qY_0)E(Y_1)\Big)
$$
$$
=
2(t-pY_1)(t-qY_1)\Big( \F(x,Y_0)-\F(Y_0,x) \Big)
+(t-px)(t-qx)\Big(\F(Y_0,Y_1)-\F(Y_1,Y_0) \Big).
$$
As above, we use the expression of the kernel to express the
coefficient of $E(x)$ as a Laurent polynomial in $x$ and $t$, and
split the right-hand side into a symmetric 
and an anti-symmetric part. After multiplying by $x$, we obtain:
$$
2\frac{(t-(1-p)x+t^2qrx^2)(t-(1-q)x+t^2pr x^2)}{pqr^2x^3}
E(x) + \frac{t(p-q)I(x)}{p^2q^2rx^2} = \hskip 5cm
$$
\beq
x(t-px)(t-qx)
\Big((t-pY_1)(t-qY_1)E(Y_0)+(t-pY_0)(t-qY_0)E(Y_1)\Big)
+(p-q)\frac {t^3B(x)}{pq}
\label{E-main}
\eeq
where $B(x)$ is given by~\Ref{Tdef-new} and
$$
I(x)= -x^3pq-t^3pqrx^3+pqr(1-2r)t^2x^4+(r+3pq)x^2t-(1+r)xt^2+t^3.
$$
As a Laurent series in $x$, the left-hand side of the above identity  has
valuation $-2$. The term involving $E(Y_0)$ and $E(Y_1)$  only
involves powers of $x$ 
smaller than or equal to $2$. But $B(x)$ is a series in $t$ with
coefficients in $\qs[x,\bx]$, containing arbitrarily large positive
and negative of $x$, and this is where the transcendence of the
solution stems from.

Note that the coefficients of $x$ and $x^2$ in $B(x)$ are especially
simple, being respectively $-t$ and $-1$.
Let  $E_i\equiv E_i(t)$ denote the coefficient of $x^i$ in
$E(x)$. Given that $E(x)=xD(x,0)$ and $D(x,y)$ is antisymmetric, we
have $E_0=E_1=0$. Let us first extract from~\Ref{E-main} the
coefficient of $x$: we obtain a relation between $E_2$ and $E_3$:
$$
E_3= \frac r t E_2 +r(q-p).
$$
Now, extracting the positive part of~\Ref{E-main} and multiplying by $x$ gives
\beq \label{Esol}
2\frac{(t-(1-p)x+t^2qrx^2)(t-(1-q)x+t^2pr x^2)}{pqr^2x^2}
E(x) = x(p-q)t^3 \frac{B^+(x)}{pq} + L(x)
\eeq
where
$B^+(x)$ is the positive part of $B(x)$
and $L(x)$ is a polynomial in $x$ (of degree $3$) with coefficients in
$\qs[p,q,t,E_2,E_4]$:
$$
L(x)=
-2\ \frac{pqrx^3t^2+tx-t^2-r(1-r)x^2t^3+x^2(r^2-pq)}{pqr^2}\ E_2
+2\ \frac{x^2t^2}{pqr^2}\ E_4 \hskip 5cm
$$
\beq
\hskip 30mm
-(p-q)\frac{(1-2r)t^3 }{pq} x^3 
+(p-q)t \frac{2+2r+t^3r}{pqr} x^2- 2(p-q) \frac{t^2}{pqr} x.
 \label{Edef-new}
\eeq

We have to determine two unknown functions $E_2$ and
$E_4$. From~\Ref{SDdef-t} and the fact that $T(x)=xS(x,0)$ and
$E(x)=xD(x,0)$, we derive that $E(t/p)=T(t/p)$ and $E(t/q)=-T(t/q)$. 
We evaluate Eq.~\Ref{Esol} at $x=t/p$ and
$x=t/q$, using the expressions of $T(t/p)$ and $T(t/q)$ given
by~\Ref{Ptp} and~\Ref{Ptq}. One thus obtains expressions of $E_2$ and
$E_4$ in terms of $\sqrt{\Delta(t/p)}, \sqrt{\Delta(t/q)}, B^+(t/p)$
and $B^+(t/q)$. They become much simpler using
$$
B^+(x)= \frac{\sqrt{\Delta(x)}}{pqrt}(1-2t\bx -pqrtx^2)
-\frac{1-2t\bx}{pqt}-2C^-(\bx). 
$$
One finds:
$$
E_2=
\frac{t^2r}{1-t^3}\left(qC^-(q/t)-pC^-(p/t)\right)
$$
and
$$
E_4=
\frac{r^2}{1-t^3}
\left((q(1-q-pt^3)C^-(q/t)-p(1-p-qt^3)C^-(p/t)\right)
-\frac{r(p-q)(2r+1)}{2t}.
$$
The expression of $D(x,0)=E(x)/x$ follows from these values,
using~\Ref{Esol},~\Ref{Edef-new} and~\Ref{TUplus}.   

\medskip
Let us now discuss the algebraic nature of $D(x,0)$ (equivalently, of
$E(x)$). If $E(x)$ were 
algebraic, the so would be all the coefficients $E_i$. In particular,
the series $L(x)$ given by~\Ref{Edef-new} and occurring in the right-hand side  of~\Ref{Esol}
would be algebraic too. By 
difference, the positive part of $B(x)$, denoted above by $B^+(x)$ would
be algebraic, and so would be  all
its coefficients. But the coefficient of $x^3$ in
$B(x)$ is
$$
4\sum_{k\ge 0} t^{3k+2} (pqr)^{k+1} \frac{(3k)!}{k!(k+1)!(k+2)!}.
$$
As $k\rightarrow \infty$, the coefficient of $t^{3k+2}$ in this series
is asymptotic to $(27pqr)^k/k^4$, up to a positive multiplicative
constant. Because of the factor  $k^{-4}$,  this cannot be the
asymptotic behaviour of the 
coefficients of an algebraic series~\cite{flajolet}, so that our
initial hypothesis is false: The series $D(x,0)$ is not algebraic. 
However, the general results on D-finite series recalled at the end of
Section 1 imply that it is D-finite.
\cqfd

\noindent
{\bf Proof of Corollary~\ref{corollary-Pi0}.} Assume all the series
$P_{i,0}$ are algebraic. By symmetry, all the series $P_{0,j}$, which
count walks ending on the $y$-axis, are algebraic too. In other words,
the coefficient of $x^i$ in $Q(x,0)$ and $Q(0,x)$ is algebraic. In
view of~\Ref{SDdef-t}, this holds for $S(x,0)$ and $D(x,0)$ as well. 

We already know, by Theorem~\ref{thm:alg-part}, that $S(x,0)$ is
algebraic. Let us work with $D(x,0)$ to obtain a contradiction.  We
thus assume that the coefficient of $x^i$ in the series $E(x)=xD(x,0)$
is algebraic. By~\Ref{Esol}, this implies that the coefficient of
$x^i$ in $B^+(x)$ is algebraic too, for all $i$. The same asymptotic
argument as above proves that this is wrong.
\cqfd

\noindent
{\bf Acknowledgements.} 
The story of this paper started when Roland Bacher
re-discovered experimentally Kreweras' result for walks ending at the
origin, and advertised his conjecture. I found it nice and
started to advertise it too. I am very grateful to Ira
Gessel who told me that this ``conjecture'' had already been proved
four times (at least), and indicated the right references. This paper
has benefited from discussions with many colleagues, in particular 
\'Eric Amar, Jean Berstel, Guy Fayolle,  Serguei Fomin, Jean-Fran\c cois
Marckert, Marni Mishna and
Nicolas Pouyanne. I thank them warmly for their interest and their
help with complex analysis, context-free languages, functional
equations, Markov chains, English and what not.



\bibliographystyle{plain}
\bibliography{qdp}

\end{document}